\documentclass[12pt]{article}

\usepackage{subcaption}

\usepackage{float}

\usepackage{multirow}

\usepackage{authorindex} 

\usepackage{tikz}

\usetikzlibrary{calc,backgrounds,arrows,matrix}

\usepackage{enumerate}

\usepackage{color}
\definecolor{lightblue}{rgb}{0,0.2,0.5}

\usepackage[colorlinks=true, urlcolor=blue,linkcolor=blue, citecolor=lightblue]{hyperref}

\usepackage[round,sort,comma,numbers]{natbib}

\bibpunct{\textcolor{lightblue}{(}}{\textcolor{lightblue}{)}}{,}{a}{}{;}

\usepackage{amssymb,amsmath}

\usepackage{graphicx}

\usepackage{graphics}
\usepackage{color} 

\DeclareMathAlphabet{\eufrak}{U}{}{}{}
\SetMathAlphabet\eufrak{normal}{U}{euf}{m}{n}
\SetMathAlphabet\eufrak{bold}{U}{euf}{b}{n}

\makeatletter
\newcommand*\rel@kern[1]{\kern#1\dimexpr\macc@kerna}
\newcommand*\widebar[1]{
  \begingroup
  \def\mathaccent##1##2{
    \rel@kern{0.8}
    \overline{\rel@kern{-0.8}\macc@nucleus\rel@kern{0.2}}
    \rel@kern{-0.2}
  }
  \macc@depth\@ne
  \let\math@bgroup\@empty \let\math@egroup\macc@set@skewchar
  \mathsurround\z@ \frozen@everymath{\mathgroup\macc@group\relax}
  \macc@set@skewchar\relax
  \let\mathaccentV\macc@nested@a
  \macc@nested@a\relax111{#1}
  \endgroup
}
\makeatother

\allowdisplaybreaks

 \def\qu{{\mathord{\mathbb Z}}}

 \def\Var{{\mathrm{{\rm Var}}}}

 \def\inte{{\mathord{\mathbb R}}}

 \def\inte{{\mathord{\mathbb N}}}

 \def\sZZ{{\rm Z\kern-.45em{}Z}}

 \def\sQQ{{\kern 0.27em \vrule height1.45ex width0.03em depth0em
           \kern-0.30em \rm Q}}
 \def\qu{{\mathchoice
         {\sQQ}
         {\sQQ}
   {\kern 0.225em \vrule height1.05ex width0.025em depth0em \kern-0.25em \rm Q}
   {\kern 0.180em \vrule height0.78ex width0.020em depth0em \kern-0.20em \rm Q}
         }}
 \def\sGG{{\kern 0.27em \vrule height1.45ex width0.03em depth0em
           \kern-0.30em \rm G}}
 \def\gg{{\mathchoice
         {\sGG}
         {\sGG}
   {\kern 0.225em \vrule height1.05ex width0.025em depth0em \kern-0.25em \rm G}
   {\kern 0.180em \vrule height0.78ex width0.020em depth0em \kern-0.20em \rm G}
         }}

 \newtheorem{prop}{Proposition}[section]

 \newtheorem{corollary}[prop]{Corollary}
 
 \newtheorem{remark}[prop]{Remark}

\numberwithin{equation}{section}

 \def\P{{\mathord{\mathbb P}}}

\newcommand{\re}{\mathrm{e}}

 \newcounter{hyp}
\usepackage{aeguill}

\usepackage{dsfont}

\newenvironment{Proofy}{\removelastskip\par\medskip
\noindent{\em Proof of Proposition} \rm}{\penalty-20\null\hfill$\square$\par\medbreak}

\newenvironment{Proof}{\removelastskip\par\medskip \noindent{\em Proof.} \rm}{\penalty-20\null\hfill$\square$\par\medbreak}

\def\bprf{\begin{Proof}}
\def\nprf{\end{Proof}}
\def\bdes{\begin{description}}
\def\ndes{\end{description}}

\newtheorem{thm}{Theorem}[section]

\def\bdef{\begin{defn}}
\def\ndef{\end{defn}}
\def\bthm{\begin{thm}}
\def\nthm{\end{thm}}
\def\bprop{\begin{prop}}
\def\nprop{\end{prop}}
\def\brmk{\begin{remark}}
\def\nrmk{\end{remark}}
\def\bexa{\begin{exa}}
\def\nexa{\end{exa}}
\def\blem{\begin{lem}}
\def\nlem{\end{lem}}
\def\bcor{\begin{cor}}
\def\ncor{\end{cor}}
\def\bexe{\begin{exe}}
\def\nexe{\end{exe}}

\newcommand{\E}{\mathbb{E}}

\newcommand{\real}{\mathbb{R}}

\def\Var{\mathop{\hbox{\rm Var}}\nolimits}

\def\og{\leavevmode\raise.3ex
     \hbox{$\scriptscriptstyle\langle\!\langle$~}}
\def\fg{\leavevmode\raise.3ex
     \hbox{~$\!\scriptscriptstyle\,\rangle\!\rangle$}~}

\title{\Huge
 Moments of Markovian growth-collapse processes 
}

 \author{Nicolas Privault
\\ 
\small
Division of Mathematical Sciences 
\\ 
\small 
School of Physical and Mathematical Sciences 
\\ 
\small
Nanyang Technological University 
\\ 
\small 
21 Nanyang Link 
\\ 
\small
Singapore 637371
}

\usepackage{hyphenat}
 
\begin{document}

\maketitle 

\vspace{-0.9cm}

\baselineskip0.6cm
 
\begin{abstract} 
We apply general moment identities for Poisson stochastic integrals with random integrands to the computation of the moments of Markovian growth-collapse processes. This extends existing formulas for mean and variance available in the literature to closed form moments expressions of all orders. In comparison with other methods based on differential equations, our approach yields polynomial expressions in the time parameter. We also treat the case of the associated embedded chain.
\end{abstract} 
 
\noindent {\bf Key words:} 
Growth-collapse processes,
Poisson shot noise,
uniform cut-off,
stochastic integrals with jumps,
moments,
cumulants.
 
\noindent
    {\em Mathematics Subject Classification (2010):}
    60J27; 
    60G55; 
    60J22. 
    
\baselineskip0.7cm

\section{Introduction}
 Markovian growth-collapse processes, see \cite{eliazar}, are
 piecewise\hyp{}deterministic Markov processes (\cite{davispdmp}),
 that grow in between random jump times at which they may
 randomly crash.
 Growth-collapse processes are used in e.g. earth sciences and physics,
 and they have also been recently applied to the study of crypto-currencies,
see \cite{frolkova}. 

\medskip 

The computation of moments of growth-collapse processes has been
the object of several approaches, see
\cite{boxma} for the use of conditional distributions
for the computation of mean and variance, and
\cite{daw} for moment expressions of all orders using
 the solution of differential equations by matrix exponentials.

\medskip 

In this paper,
we apply general moment identities
 written as sums over partitions 
for Poisson stochastic integrals with
random integrands, see \cite{momentpoi,prinv,priinvcr,prob},
 to the computation of the moments of growth-collapse processes.
In particular, we obtain closed-form moment expressions which are polynomial
in time in the case of uniformly distributed cut-off rates. 

\medskip 

Let $(N_t)_{t\in \real_+}$ denote a standard Poisson process
with intensity $\lambda$ on $\real_+$, and
consider a process $(Y_t)_{t\in \real_+}$ of the form 
$$ 
 Y_t  = \int_0^t f(N_{s^-}) dN_s,
 \qquad t\geq 0. 
$$
 As the left limit $f(N_{s^-})$ is predictable with respect to the
 filtration $({\cal F}_s)_{s\in \real_+}$ generated by $(N_s)_{s\in \real_+}$
 the mean of $Y_t$
 can be computed from the smoothing lemma, see e.g. Theorem~9.2.1 in \cite{Bh},
 as 
$$ 
 \E [ Y_t ] = 
 \lambda \E \left[ \int_0^t f(N_{s^-}) ds \right]
 =
 \lambda \int_0^t \E \left[ f(N_{s^-}) \right] ds, 
 \qquad t\geq 0. 
$$
 This calculation does not apply however to the process
 $f(N_s)$ which is only adapted and not predictable with respect to 
 the filtration $({\cal F}_s)_{s\in \real_+}$.
 In this case we may apply the 
 Slivnyak-Mecke formula, see \cite{slivnyak}, \cite{jmecke}, 
 to obtain 
$$ 
 \E [ Y_t ] = 
 \lambda \E \left[ \int_0^t \varepsilon^+_s f(N_s) ds \right]
 =
 \lambda \E \left[ \int_0^t f(1+N_s) ds \right], 
 \qquad t\geq 0, 
$$
 where $\varepsilon^+_s$ denotes the operator that adds one jump at the
 location $s\geq 0$ to the Poisson process path.

 \medskip

 In order to compute higher order moments we will
 apply a nonlinear extension of the 
 Slivnyak-Mecke identity, see Proposition~\ref{fldsf} below,
 which allows us to express the moments of Poisson stochastic integrals
 as a sum of multiple integrals with respect to the intensity of the Poisson
 process over partitions.
 In Section~\ref{s3} we consider the computation of moments of
 jump processes of the form
$$ 
  Y_t = \sum_{k=1}^{N_t} g(T_k,k,N_t)
 = \int_0^t g(s,N_s,N_t) dN_s, 
$$
 where $(T_k)_{k\geq 1}$ denotes the sequence of jump times
 of the Poisson process $(N_t)_{t\in \real_+}$,
 see Proposition~\ref{p3.1} and its Corollary~\ref{p3.2}.
 Those identities are then specialized
 in Section~\ref{s4} to the case of uniform cut-off
 distributions, for processes of the form
 $$
 Y_t  = \sum_{k=1}^{N_t} f_k (T_k) (1-U_k) \prod_{l=k+1}^{N_t} U_l,
\qquad t\in \real_+, 
$$
 where $(U_k)_{k\geq 1}$ is an i.i.d. uniform sequence on $[0,1]$,
 independent of the standard Poisson process $(N_t)_{t\in \real_+}$,
 see Corollary~\ref{p4.1}.  

 \medskip

 In particular, in Proposition~\ref{p5.1} we obtain the 
 closed form polynomial expression 
 \begin{equation}
   \label{expr} 
  \E [ (X_t)^n ] = \frac{(n+1)!}{\lambda^n} \sum_{k=0}^n (-1)^k (k+1)^{n-1}{n \choose k}
  \re^{-k\lambda t/(k+1)},
  \qquad t\in \real_+, 
\end{equation} 
 for the moments of all orders $n\geq 0$ of the
 growth-collapse process 
 $$
X_t = t - \sum_{k=1}^{N_t} T_k (1-U_k) \prod_{l=k+1}^{N_t} U_l,
 \qquad t\in \real_+, 
$$
 where $(U_k)_{k\geq 1}$ is an i.i.d. uniform sequence on $[0,1]$.
 This result extends Theorems~4 and 5
 as well as Corollary~1 of \cite{boxma} from
 mean and variance to higher moments of all orders,
 and provides a closed form alternative to
 Corollary~4 in \cite{daw} which uses 
 matrix exponentials.
 The expression \eqref{expr} immediately yields the asymptotic moments 
$$
  \lim_{t\to \infty} 
  \E [ (X_t)^n ] = \frac{(n+1)!}{\lambda^n}, 
  \qquad n\geq 1, 
$$
 which recover  
 the gamma stationary distribution of $(X_t)_{t\in \real_+}$
 with shape parameter $2$, see Theorem~3 in \cite{boxma}. 

 \medskip

 Finally, in Section~\ref{s6} we show that our approach 
 can also applied to discrete-time embedded processes of
 the form 
$$ 
  Y(m) = \sum_{k=1}^m g(T_k,k,m)
  = \int_0^{T_m} g(s,N_s,m) dN_s,
  \qquad m \geq 1, 
$$
  see Corollaries~\ref{p6.2}-\ref{p6.3}, 
  and to the embedded growth-collapse chain 
$$
 X(m) = T_m - \sum_{k=1}^m T_k (1-U_k) \prod_{l=k+1}^m U_l,
 \qquad m \geq 1,  
$$
 where $(U_k)_{k\geq 1}$ is an i.i.d. uniform sequence on $[0,1]$, 
 see Corollaries~\ref{p6.4}-\ref{p6.5}. 
 This recovers Theorem~7 stated for mean and variance
 in \cite{boxma},
 and provides moment expressions of all orders.

 \medskip  
 
 We proceed as follows.
In Section~\ref{s2} we review the derivation
of moment identities for stochastic integrals using sums over partitions,
and in Section~\ref{s3} we apply them to the moments of
 jump processes driven by a Poisson process.
 Those expressions are then specialized as closed form polynomial
 identities in Section~\ref{s4} in the case of uniform cut-off
 distributions. 
 The moments of growth-collapse processes 
 are considered in Section~\ref{s5}, and the case
 of embedded chains is treated in Section~\ref{s6}. 

\section{Moment identities for Poisson stochastic integrals} 
\label{s2}
In this section we review the computation of moments of
Poisson stochastic integrals with random integrands using sums over partitions,
see Proposition~3.1 in \cite{momentpoi}. 
Consider a Poisson process $(N_t)_{t\in \real_+}$ constructed as
$N_t = \omega ([0,t])$, where $\omega (ds)$ is a Poisson random
measure of intensity $\lambda (ds)$, 
with sequence $(T_k)_{k\geq 1}$ of jump times.
 For any $s_1,\ldots ,s_k\in \real_+$, we let
 $\epsilon^+_{s_1} \cdots \epsilon^+_{s_k}$ denote  
 the operator 
$$ 
 ( \epsilon^+_{s_1} \cdots \epsilon^+_{s_k} F)(\omega  )=F(\omega  \cup\{s_1, \ldots, s_k\})
 $$
 acting on random variables $F$ by addition of points
 at locations $s_1, \ldots, s_k$ to the point process $\omega  (dx)$.
 For example, if $F$ takes the form
 $F = f(N_{t_1},\ldots , N_{t_n})$, then we have 
$$ 
 \epsilon^+_{s_1} \cdots \epsilon^+_{s_k} F
 =
 f(N_{t_1} + \# \{ k \ : \ s_k \leq t_1\}
 ,\ldots , N_{t_n} + \# \{ k \ : \ s_k \leq t_n\}). 
$$
The following moment identity, see Proposition~3.1 in \cite{momentpoi} 
and Theorem~1 in \cite{prob},
 uses sums over partitions $\{\pi_1,\ldots , \pi_k \}$
 of $\{1,\ldots , n\}$, and applies to random integrands 
 $u : \real_+ \times \Omega \longrightarrow \real$.
 \begin{prop}
   \label{fldsf}
 Let $(u_s(\omega  ))_{s\in \real_+}$ denote a stochastic process
 indexed by $s\in \real_+$.
 For any $n\geq 1$, we have 
\medskip 
\begin{align} 
\label{dkjlddss} 
& \E \left[ 
  \left(
  \int_0^T u_s dN_s 
  \right)^n
 \right] 
\\ 
\nonumber 
 & = 
\sum_{k=1}^n 
\sum_{ \pi_1 \cup \cdots \cup \pi_k = \{ 1 , \ldots , n\} } 
 \E \left[ 
 \int_0^t \cdots \int_0^t 
 \epsilon^+_{s_1} \cdots \epsilon^+_{s_k} 
 \big( 
 u^{|\pi_1|} (s_1, \omega ) 
 \cdots 
 u^{|\pi_k|} ( s_k, \omega ) 
 \big) 
 \lambda ( d s_1 ) \cdots \lambda ( d s_k ) 
 \right],  
\end{align} 
where the power $|\pi_i|$ denotes the cardinality of the subset $\pi_i$
and the above sum runs over all partitions 
$\pi_1,\ldots , \pi_k$ of $\{ 1 , \ldots , n \}$.
\end{prop} 
 In the sequel we will frequently use the equivalent
 combinatorial expressions
\begin{eqnarray*} 
 \sum_{k=1}^n 
 \sum_{\pi_1 \cup \ldots \cup \pi_k = \{1,\ldots , n \} } 
 f_k(|\pi_1|, \ldots , |\pi_k|)
   & = & 
   \sum_{k=1}^n 
  \frac{n!}{k!}
  \sum_{p_1+\cdots + p_k=n \atop p_1,\ldots, p_k \geq 1} 
 \frac{f_k(p_1, \ldots , p_k) }{p_1! \cdots p_k!}
  \\
   & = & 
  \sum_{k=1}^n 
 \frac{n!}{k!} 
  \sum_{q_0=0< q_1 < \cdots < q_k = n } 
  \frac{f_k(q_1-q_0, \ldots , q_k-q_{k-1})}{
    (q_1-q_0)! \cdots (q_k-q_{k-1})!    }
\end{eqnarray*} 
 for $f_k$ a function on $\inte^k$, $k=1,\ldots , n$.  
 In particular,
 for $x_1,\ldots , x_n\in \real$ and
 $f_k(p_1,\ldots , p_k) = x_{p_1}\cdots x_{p_k}$
 this yields the Bell polynomial of order $n\geq 1$ as 
\begin{eqnarray} 
\nonumber 
 B_n ( x_1 , \ldots , x_n ) 
 & = & 
 \sum_{k=1}^n
 \frac{n!}{k!} 
 \sum_{p_1+\cdots + p_k=n \atop 
 p_1\geq 1,\ldots ,p_k \geq 1 
 } 
 \frac{x_{p_1} \cdots x_{p_k}}{p_1! \cdots p_k!} 
\\ 
\nonumber 
 & = & 
 \sum_{k=0}^n 
 ~ 
 \sum_{\pi_1\cup \cdots \cup \pi_k = \{1,\ldots ,n \} } 
 x_{|\pi_1|} 
 \cdots 
 x_{|\pi_k |} 
. 
\end{eqnarray} 
We will also use the relation 
$\E [ X^n ] 
  = 
 B_n ( 
 \kappa_X^{(1)} 
 ,
 \ldots 
 , 
 \kappa_X^{(n)} 
 )$ 
 between the moments
 $\E [ X^n]$ and the cumulants
 $\kappa_X^{(n)}$ of a random variable $X$,
 and the inversion relation 
 \begin{equation}
   \label{inv} 
   \kappa^X_n = \sum_{k=1}^n 
 (k-1)! (-1)^{k-1} 
 \sum_{\pi_1 \cup \cdots \cup \pi_k = \{ 1,\ldots ,n \} } 
 \E \big[ X^{|\pi_1|} \big] \cdots \E \big[ X^{|\pi_k|} \big],
 \qquad n \geq 1, 
\end{equation} 
 see Theorem~1 of \cite{elukacs}, or \cite{leonov}. 
\subsubsection*{Shot noise processes}
  Before moving to the setting of Markovian growth-collapse processes,
we use the case of Poisson shot noise processes as an illustration for
the result of Proposition~\ref{fldsf}.
 Consider a shot noise process $(S_t)_{t\in \real_+}$ of the form
$$
S_t = \sum_{k=1}^{N_t} J_k g(T_k,t)
= \sum_{k=1}^{N_t} J_{N_{T_k}} g( T_k , t)
= \int_0^t J_{N_s} g( s , t) dN_s,
\qquad t\in \real_+, 
$$
where $(J_k)_{k\geq 0}$ is a sequence of i.i.d. random variables
admitting moments of all orders,
and $g(\cdot , \cdot )$ is a sufficiently integrable
deterministic function. 
The next proposition provides a closed form expression for
the moments of shot noise processes using standard Bell polynomials,
see also Corollary~2 in \cite{daw} for another expression using matrix
exponentials
in case $\lambda (ds) = \lambda ds$ for some rate $\lambda >0$,
 and $g(s,t) = \re^{-\beta (t-s)}$ for some $\beta >0$.  
\begin{prop}
  \label{p1}
  For any $n\geq 1$, we have 
$$ 
 \E [ S_t^n ] 
 = B_n \left(
\E \big[
   J_1
   \big]
 \int_0^t g (s,t) ds  
 ,
 \ldots , 
 \E \big[
   J_1^n
   \big]
 \int_0^t g^n(s,t) ds  
 \right)
 ,
 $$
 where $B_n$ is the Bell polynomial of order $n\geq 1$.
\end{prop}
\begin{Proof}
  Taking
  $u_s (\omega ) : = J_{N_s} g(s,t)$, 
  by \eqref{dkjlddss} we have
\begin{align*} 
\nonumber 
& \E [ S_t^n ] 
 = 
\E \left[ \left( \sum_{k=1}^{N_t} J_{N_{T_k}} g(T_k,t) \right)^n \right]
\\ 
\nonumber 
 & = 
\E \left[ \left( \sum_{k=1}^{N_t} u_{T_k} (\omega ) \right)^n \right]
\\ 
\nonumber 
 & = 
\sum_{l=1}^n 
\sum_{ \pi_1 \cup \cdots \cup \pi_l = \{ 1 , \ldots , n\} } 
\E \left[
  \int_0^t \cdots \int_0^t 
 \epsilon^+_{s_1} \cdots \epsilon^+_{s_l} 
 \big( 
 u^{|\pi_1|} (s_1, \omega ) 
 \cdots 
 u^{|\pi_l|} ( s_l, \omega ) 
 \big) 
 \lambda ( d s_1 ) \cdots \lambda ( d s_l )
 \right]
\\ 
\nonumber 
 & = 
\sum_{l=1}^n 
\sum_{ \pi_1 \cup \cdots \cup \pi_l = \{ 1 , \ldots , n\} } 
  \int_0^t \cdots \int_0^t 
 g^{|\pi_1|}(s_1,t) 
 \cdots 
 g^{|\pi_l|} (s_l,t)
 \E \big[
 \epsilon^+_{s_1} \cdots \epsilon^+_{s_l} 
 \big( 
 J_{N_{s_1}}^{|\pi_1|} \cdots J_{N_{s_l}}^{|\pi_l|} 
 \big)
 \big] 
 \lambda ( d s_1 ) \cdots \lambda ( d s_l )
  \\ 
\nonumber 
 & = 
\sum_{l=1}^n 
\sum_{ \pi_1 \cup \cdots \cup \pi_l = \{ 1 , \ldots , n\} } 
 \int_0^t \cdots \int_0^t 
 g^{|\pi_1|}(s_1,t) 
 \cdots 
 g^{|\pi_l|} (s_l,t)
 \E \big[
   J_{N_{s_1}}^{|\pi_1|}
   \big]
 \E \big[
   J_{1+N_{s_1}}^{|\pi_1|}
   \big]
 \cdots
 \E \big[
   J_{l-1+N_{s_l}}^{|\pi_l|} 
   \big]
 \lambda ( d s_1 ) \cdots \lambda ( d s_l )
  \\ 
\nonumber 
 & = 
\sum_{l=1}^n 
\sum_{ \pi_1 \cup \cdots \cup \pi_l = \{ 1 , \ldots , n\} } 
 \E \big[
   J_1^{|\pi_1|}
   \big]
 \cdots
 \E \big[
   J_1^{|\pi_l|} 
   \big]
 \int_0^t \cdots \int_0^t 
 g^{|\pi_1|}(s_1,t) 
 \cdots 
 g^{|\pi_l|} (s_l,t)
 \lambda ( d s_1 ) \cdots \lambda ( d s_l ). 
\end{align*} 
\end{Proof}
 We note that Proposition~\ref{p1} is consistent with the L\'evy-Khintchine
 formula for compound Poisson processes, as the Fa\`a di Bruno formula,
 see e.g. \S~2 of \cite{elukacs}, yields
 \begin{eqnarray*}
 \E [ \re^{\alpha S_t} ] 
 & = & 
 \sum_{n\geq 0}\frac{\alpha^n}{n!}\E [ S_t^n ] 
 \\
 & = & 
 \sum_{n\geq 0}\frac{\alpha^n}{n!}
 B_n \left(
\E \big[
   J_1
   \big]
 \int_0^t g (s,t) ds  
 ,
 \ldots , 
 \E \big[
   J_1^n
   \big]
 \int_0^t g^n(s,t) ds  
 \right)
 \\
 & = & 
 \exp \left(
 \sum_{n\geq 1}
 \frac{\alpha^n}{n!}
 \E \big[
   J_1^n
   \big]
 \int_0^t g^n(s,t) ds  
 \right)
 \\
 & = & 
 \exp \left(
 \int_0^t \big(
 \re^{\alpha g (s,t) J_1}
 - 1
 \big)
 ds  
 \right), 
\end{eqnarray*} 
 which recovers the cumulants of $S_t$ from the moments of $J_1$ as 
$$
 \kappa_{S_t}^{(n)} =
 \E \big[
   J_1^n \big]
 \int_0^t g^n (s,t) ds,
 \qquad n \geq 1.
$$ 
\section{Moments of jump processes} 
\label{s3}
\noindent
From now on we assume that $(N_t)_{t\in \real_+}$ is a standard Poisson process
with intensity $\lambda >0$, and in this section we consider jump processes
built as the anticipating Poisson integrals 
\begin{equation} 
  \label{dfjkl}
  Y_t = \sum_{k=1}^{N_t} g(T_k,k,N_t)
 = \int_0^t g(s,N_s,N_t) dN_s, \qquad t\in \real_+. 
\end{equation} 
  
\begin{prop}
\label{p3.1}
 Let $(Y_t)_{t\in \real_+}$ be defined as in \eqref{dfjkl}.
 For all $n\geq 1$, we have
$$ 
 \E [ (Y_t)^n] 
=
 \sum_{k=1}^n 
 ~ 
 \lambda^k
 \sum_{\pi_1 \cup \ldots \cup \pi_k = \{1,\ldots , n \} } 
  \int_0^t \cdots \int_0^t
 \E \left[ 
    \prod_{l=1}^k g^{|\pi_l|} (s_l,l+N_{s_l},k+N_t)
   \right] 
 ds_1 \cdots ds_n
. 
$$
\end{prop}
\begin{Proof}
 We have 
\begin{eqnarray*} 
   \E [ (Y_t)^n] & = &  
  \E \left[
 \left(   \int_0^t g(s,N_s,N_t) dN_s
    \right)^n  \right]
      \\
& = & 
 \sum_{k=1}^n 
 ~ 
 \lambda^k
 \sum_{\pi_1 \cup \ldots \cup \pi_k = \{1,\ldots , n \} } 
  \int_0^t \cdots \int_0^t
 \E \left[ 
  \epsilon^+_{s_1} \cdots \epsilon^+_{s_n}
 \prod_{l=1}^k
 g^{|\pi_l|} (s_l,N_{s_l},N_t) 
 \right] 
 ds_1 \cdots ds_n
      \\
& = & 
 \sum_{k=1}^n 
 ~ 
 \lambda^k
 \sum_{\pi_1 \cup \ldots \cup \pi_k = \{1,\ldots , n \} } 
  \int_0^t \cdots \int_0^t
 \E \left[ 
    \prod_{l=1}^k g^{p_l} (s_l,l+N_{s_l},k+N_t)
   \right] 
 ds_1 \cdots ds_n, 
\end{eqnarray*} 
where the sum runs over all partitions 
$\pi_1,\ldots , \pi_k$ of $\{ 1 , \ldots , n \}$.
\end{Proof}
Next, we specialize Proposition~\ref{p3.1} to the case where $g(s,k,n)$ takes the
form
$$
g(s,k,n) = f_k(s) W_k^n, 
$$
 where 
$$
 W_k^n = (1-Z_k) \prod_{l=k+1}^n Z_l,
 \qquad 1 \leq k \leq n, 
$$
 and $(Z_k)_{k \geq 1}$ is an i.i.d. random sequence independent
of the Poisson process $(N_t)_{t\in \real_+}$, with moment sequence
$m_n = \E [ Z^n]$, $n\geq 0$, i.e. we have 
\begin{equation}
\label{dfjkl2}
Y_t  = \sum_{k=1}^{N_t} f_k (T_k) (1-Z_k) \prod_{l=k+1}^{N_t} Z_l,
\qquad t\in \real_+. 
\end{equation} 
  
\begin{corollary}
  \label{p3.2}
  Let $(Y_t)_{t\in \real_+}$ be defined as in \eqref{dfjkl2}.
  For all $n\geq 1$ we have 
\begin{eqnarray*} 
\lefteqn{ 
 \E [ (Y_t)^n]
=  n!
 \re^{\lambda t (m_n-1)} 
  \sum_{k=1}^n 
  \lambda^k
  \hskip-0.4cm
  \sum_{q_0=0< q_1 < \cdots < q_k = n } 
}
  \\
  & & \qquad \qquad \quad
  \int_0^t \int_0^{s_k} \cdots \int_0^{s_2}
 \prod_{l=1}^k
 \left(
 \frac{f_l^{q_l-q_{l-1}}(s_l)}{
   (q_l-q_{l-1})!} 
 C_{q_{l-1},q_l-q_{l-1}}
 \re^{\lambda s_l ( m_{q_{l-1}} - m_{q_l} ) } 
 \right)
 ds_1 \cdots ds_k
, 
\end{eqnarray*} 
 $t\in \real_+$, where 
$$
 C_{p,q}
 := 
 \E [ ( 1 - Z )^p Z^q]
 = 
 \sum_{k=0}^p {p \choose k} (-1)^k
 m_{q+k}, \qquad p,q \geq 0. 
 $$
\end{corollary}
\begin{Proof}
  By Proposition~\ref{p3.1},
  for all $n\geq 1$ we have
\begin{eqnarray*} 
  \lefteqn{
    \! \! 
    \E [ (Y_t)^n] = 
  n!
  \sum_{k=1}^n 
  \lambda^k
 \sum_{p_1+\cdots + p_k=n \atop p_1,\ldots, p_k \geq 1} 
}
\\
& &
\qquad  \quad 
\int_0^t \int_0^{s_k} \cdots \int_0^{s_2}
 \frac{f^{p_1}(s_1) \cdots f^{p_k}(s_k)}{p_1!\cdots p_k!}
 \E \left[ 
 \epsilon^+_{s_1} \cdots \epsilon^+_{s_n}
 \big( 
  W_{N_{s_1}}^{N_t} 
 \big)^{p_1} 
 \cdots
 \big(
  W_{N_{s_k}}^{N_t} 
 \big)^{p_k} 
 \right] 
 ds_1 \cdots ds_k
. 
\end{eqnarray*} 
 For $p_1+\cdots + p_k=n$ and $0\leq s_1<\cdots < s_k\leq s_{k+1}:=t$, we have 
\begin{eqnarray*} 
  \lefteqn{
    \! \! \! \! \! \! \! \! \! \! \! \! \! \! \! \! \! \! \! \! 
    \E \left[
    \epsilon^+_{s_1} \cdots
    \epsilon^+_{s_k}
    \Big( 
 \big(   W_{N_{s_1}}^{N_t} 
 \big)^{p_1}
 \cdots 
 \big( W_{N_{s_k}}^{N_t} 
  \big)^{p_k}
  \Big)
  \right]
    =
    \E \left[
        \big(
    W_{1+N_{s_1}}^{k+N_t} 
    \big)^{p_1}
    \cdots 
    \big(
    W_{k+N_{s_k}}^{k+N_t} 
  \big)^{p_k}
    \right]
  }
  \\
& = &
  \E \left[
    \left(
    (1-Z_{1 + N_{s_1}}) \prod_{l=2+N_{s_1}}^{k+N_t} Z_l
    \right)^{p_1}
    \cdots 
    \left(
    (1-Z_{k+N_{s_k}}) \prod_{l=k+1+N_{s_k}}^{k+N_t} Z_l
    \right)^{p_k}
    \right]
 \\
& = &
 \E \left[
 \prod_{l=1}^k
 \left(
 ( 1 - Z_{l+N_{s_l}} )^{p_l}
 ( Z_{l+N_{s_l}} )^{p_1+\cdots + p_{l-1}}
 \prod_{p=1+N_{s_l}}^{N_{s_{l+1}}}
 ( Z_{l+p} )^{p_1+\cdots + p_l}
 \right)
 \right]
\\
& = &
 \E \left[
 \prod_{l=1}^k
 \left(
 C_{p_1+\cdots + p_{l-1},p_l}
 ( m_{p_1+\cdots + p_l} )^{N_{s_{l+1}} - N_{s_l}} 
 \right)
 \right]
\\
& = &
 \prod_{l=1}^k
 \left(
 C_{p_1+\cdots + p_{l-1},p_l}
 \re^{\lambda ( s_{l+1} - s_l )
 ( m_{p_1+\cdots + p_l} -1) }
 \right)
\\
& = &
  \re^{\lambda t ( m_n -1) }
 \prod_{l=1}^k
 \big(
 C_{p_1+\cdots + p_{l-1},p_l}
   \re^{\lambda s_l ( m_{p_1+\cdots + p_{l-1}} - m_{p_1+\cdots + p_l} ) }
\big).
\end{eqnarray*} 
\end{Proof}
\section{Uniform cut-offs} 
\label{s4}
 In this section we assume that $Z_k$ is uniform on
 $[0,1]$, $k\geq 1$. 
  In this case we have $m_n=1/(n+1)$ and
  $C_{p,q}$ is given by the beta function as
  $C_{p,q} = p!q!/(p+q+1)!$, hence 
$$ 
 \prod_{l=1}^k
 C_{p_1+\cdots + p_{l-1},p_l}
 = 
 \prod_{l=1}^k
\frac{(p_1+\cdots + p_{l-1})!p_l!}{
(p_1+\cdots + p_l+1)! 
}
=
 \frac{1}{n!} \prod_{l=1}^k
 \frac{p_l!
 }{p_1+\cdots + p_l+1}, 
$$ 
 which yields the next consequence of Corollary~\ref{p3.2}. 
\begin{corollary}
   \label{p4.1}
   Let $(Y_t)_{t\in \real_+}$ be defined as in \eqref{dfjkl2},
  where $(Z_k)_{k\geq 1}$ is an i.i.d. uniform sequence on $[0,1]$. 
  For any $n\geq 1$, we have 
  \begin{equation}
    \label{eq4.1} 
 \E [ (Y_t)^n] 
 = \re^{ - n \lambda t / (n+1) }
  \sum_{k=1}^n 
 \lambda^k
 \sum_{q_0=0< q_1 < \cdots < q_k = n } 
 \int_0^t \int_0^{s_k} \cdots  \int_0^{s_2}
 \prod_{l=1}^k
  \frac{f^{q_l-q_{l-1}}_l(s_l) \re^{ \lambda s_l (1/q_{l-1}-1/q_l) }}{1+q_l}
  ds_1 \cdots ds_k, 
\end{equation} 
  $t\in \real_+$.
 \end{corollary}
 When $f(s)=s$, $s\in \real_+$, the relation 
\begin{equation} 
  \label{rec}
  \int_0^{s_2} s_1^n
\re^{ \alpha s_1 } ds_1
=
1 - \re^{\alpha {s_2}} \sum_{k=0}^n (-1)^k \frac{(\alpha {s_2})^k}{k!} ,
\qquad
s_2\in \real_+, \quad n \geq 0, 
\end{equation} 
 can be used compute the integrals in \eqref{eq4.1} by induction. 
 \subsubsection*{First moment}
 For $n=1$ we have
\begin{equation}
\label{m1x} 
 \E [ Y_t ] = \frac{\lambda }{2} \re^{ - \lambda t / 2}  \int_0^t s_1 \re^{  \lambda s_1 / 2} ds_1 
 = t - \frac{2}{\lambda} ( 1 - \re^{-\lambda t/2} ), 
\end{equation} 
which is consistent with Theorem~4 in \cite{boxma}, with a shorter proof,
see Figure~\ref{fig1}. 
\subsubsection*{Second moment}
\begin{figure}[H]
  \centering
 \begin{subfigure}[b]{0.49\textwidth}
    \includegraphics[width=1\linewidth, height=5cm]{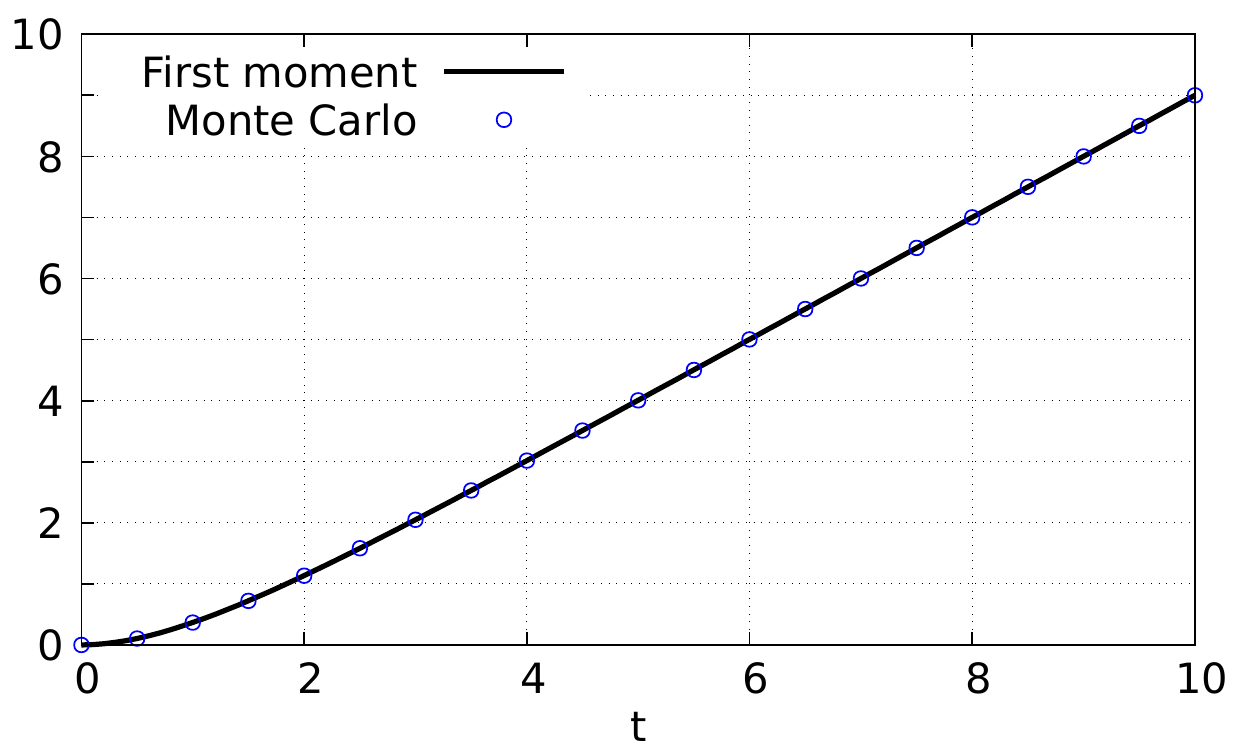}
    \caption{First cumulant $\kappa^{(1)} (t)$.} 
 \end{subfigure}
 \begin{subfigure}[b]{0.49\textwidth}
    \includegraphics[width=1\linewidth, height=5cm]{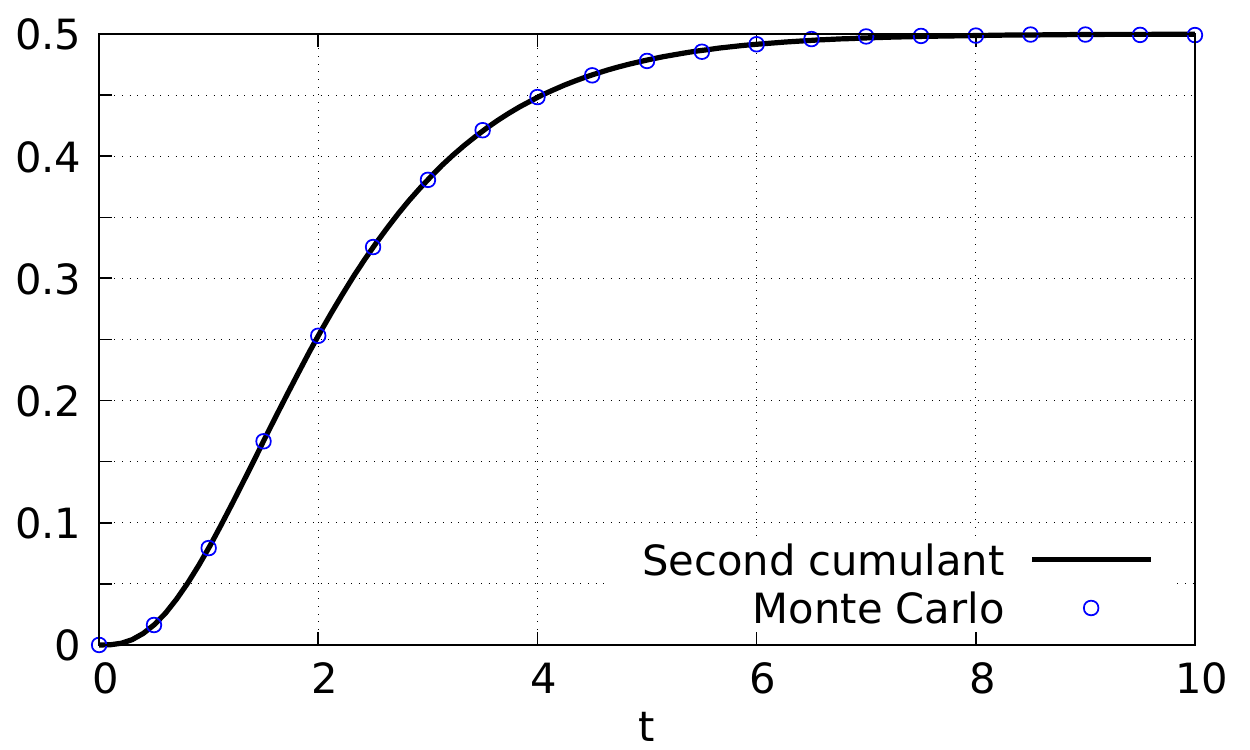} 
    \caption{Second cumulant $\kappa^{(2)} (t)$.} 
 \end{subfigure}
  \caption{Mean and variance given by \eqref{m1x} and \eqref{m2x}.} 
\label{fig1} 
\end{figure}

  For $n=2$ we have
\begin{eqnarray} 
  \nonumber
  \E [ (Y_t)^2] 
 & = &
 \frac{\lambda}{3}
 \re^{ - 2\lambda t/3}  \int_0^t 
 s_1^2
 \re^{ 2\lambda s_1 /3}  
 ds_1
 +
 \frac{\lambda^2}{6}
 \re^{ - 2 \lambda t / 3 }
  \int_0^t 
  s_2 \re^{ \lambda s_2 /6 }
  \int_0^{s_2}
  s_1 \re^{ \lambda s_1 /2} ds_1 ds_2 
  \\
  \label{m2x}
  & = &  
  \frac{18}{\lambda^2} 
  \re^{ - 2 \lambda t / 3 }
  + \frac{4}{\lambda} \re^{-\lambda t/2} t
  - \frac{24}{\lambda^2} \re^{-\lambda t/2}
 + t^2 + \frac{6}{\lambda^2} 
 - \frac{4}{\lambda} t, 
\end{eqnarray} 
 hence  
$$ 
 \kappa^{(2)}(t) = \Var [Y_t]
 = \E [ Y_t^2] - ( \E [Y_t ] )^2
 = \frac{2}{\lambda^2} \big(
 9 \re^{ - 2 \lambda t / 3 } 
 - 2 \re^{-\lambda t}
 - 8 \re^{-\lambda t/2}
 + 1 \big), 
$$ 
 which recovers Theorem~5 in \cite{boxma}
 with a shorter proof, see Figure~\ref{fig1}.
 Figures~\ref{fig1} to \ref{fig7}
 are plotted with $10$ million Monte Carlo samples
 and $\lambda =2$. 

\noindent
The subsequent integrals for higher order moments can be evaluated using Mathematica or 
based on the recurrence relation \eqref{rec}. 
\subsubsection*{Third moment}
  For $n=3$ we have
  \begin{align}
    \nonumber 
 &    \E [ (Y_t)^3] 
  = 
 \frac{\lambda}{4}
 \re^{ - 3\lambda t/4}  \int_0^t 
 s_1^3 \re^{ 3\lambda s_1 /4} ds_1
 + \frac{\lambda^2 }{8} 
  \re^{ - 3 \lambda t / 4 }
  \int_0^t 
  s_2^2 \re^{ \lambda s_2 /4 }
  \int_0^{s_2}
  s_1 \re^{ \lambda s_1 /2} ds_1 ds_2 
  \\
    \nonumber 
  &  \quad 
 + \frac{\lambda^2 }{12} 
  \re^{ - 3 \lambda t / 4 }
  \int_0^t 
  s_2 \re^{ \lambda s_2 /12}
  \int_0^{s_2}
  s_1^2 \re^{ 2\lambda s_1 /3 } ds_1 ds_2 
  \\
    \nonumber 
  &  \quad 
 +  \frac{\lambda^3 }{24} 
  \re^{ - 3 \lambda t / 4 }
  \int_0^t 
  s_3 \re^{ \lambda s_3 /12 }
  \int_0^{s_3}
  s_2 \re^{ \lambda s_2 /6 }
  \int_0^{s_2}
  s_1 \re^{ \lambda s_1 /2} ds_1 ds_2 ds_3
  \\
  \label{m3x}
  & = 
  \frac{\re^{-3 \lambda t/  4}}{\lambda^3}
  \big(384 + 54 \re^{\lambda t/12} (\lambda t -12) + 
   6 \re^{\lambda t/4} (48 + \lambda t ( \lambda t -12 )) + 
   \re^{3 \lambda t/ 4} (\lambda t (18 + \lambda t ( \lambda t -6 ))-24)
   \big),
\end{align}
 see Figure~\ref{fig2} below. 
\subsubsection*{Fourth moment}
      For $n=4$ we have
\begin{align} 
  \nonumber
  & 
  \E [ (Y_t)^4] =
  \frac{\lambda}{5} \re^{ - 4\lambda t/5}  
 \int_0^t s_1^4 \re^{4\lambda s_1/5} ds_1
 + 
 \frac{\lambda^2}{10} \re^{ - 4\lambda t/5}  
 \int_0^t s_2^3 \re^{3\lambda s_2/10} \int_0^{s_2} 
 s_1 \re^{\lambda s_1/2} ds_1 ds_2
\\
\nonumber
& \quad
 + \frac{\lambda^2}{20} \re^{ - 4\lambda t/5}  
 \int_0^t s_2 \re^{\lambda s_2/20} \int_0^{s_2} 
 s_1^3 \re^{3\lambda s_1/4} 
 ds_1 ds_2
+ 
 \frac{\lambda^2}{15} \re^{ - 4\lambda t/5}  
 \int_0^t s_2^2\re^{2\lambda s_2/15} \int_0^{s_2} 
 s_1^2 \re^{2\lambda s_1/3}
 ds_1 ds_2
 \\
\nonumber
   & \quad
 + \frac{\lambda^3}{30} \re^{ - 4\lambda t/5}  
 \int_0^t s_3^2 \re^{2\lambda s_3/15}
 \int_0^{s_3} s_2 \re^{\lambda s_2/6}
 \int_0^{s_2} s_1 \re^{\lambda s_1/2}
 ds_1 ds_2 ds_3 
 \\
\nonumber
   & \quad 
 + \frac{\lambda^3}{40} \re^{ - 4\lambda t/5}
 \int_0^t s_3 \re^{\lambda s_3/20}
 \int_0^{s_3} s_2^2 \re^{\lambda s_2/4}
 \int_0^{s_2} s_1 \re^{\lambda s_1/2}
 ds_1 ds_2 ds_3 
 \\
\nonumber
   & \quad 
  + \frac{\lambda^3}{60} \re^{ - 4\lambda t/5}  
  \int_0^t s_3 \re^{\lambda s_1/20}
  \int_0^{s_3} s_2 \re^{\lambda s_1/12}
  \int_0^{s_2} s_1^2 \re^{2\lambda s_1/3}
 ds_1 ds_2 ds_3 
 \\
\nonumber
   & \quad  
 + \frac{\lambda^4}{120} \re^{ - 4\lambda t/5}
 \int_0^t s_4 \re^{\lambda s_4/20}
 \int_0^{s_4}  s_3 \re^{\lambda s_3/12}
 \int_0^{s_3} s_2\re^{\lambda s_2/6}
 \int_0^{s_2} s_1 \re^{\lambda s_1/2}
 ds_1 ds_2 ds_3 ds_4
 \\
 \label{m4x}
 & = 
\frac{\re^{-4 \lambda t/  5}}{\lambda^4} \big(15000 + 1536 \re^{\lambda t/20} ( \lambda t -20) + 
108 \re^{2 \lambda t/15} (180 + \lambda t ( \lambda t -24))
\\
\nonumber
 & \quad + 
   8 \re^{3 \lambda t/  10} ( \lambda t (144 + \lambda t ( \lambda t -18)) -480) + 
   \re^{4 \lambda t/ 5} (120 + 
      \lambda t (\lambda t (36 + \lambda t ( \lambda t -8)) -96))\big), 
\end{align} 
 see Figure~\ref{fig3} below. 
\section{Moments of growth-collapse processes}
\label{s5}
\begin{figure}[H]
  \centering
 \begin{subfigure}[b]{0.49\textwidth}
    \includegraphics[width=1\linewidth, height=5cm]{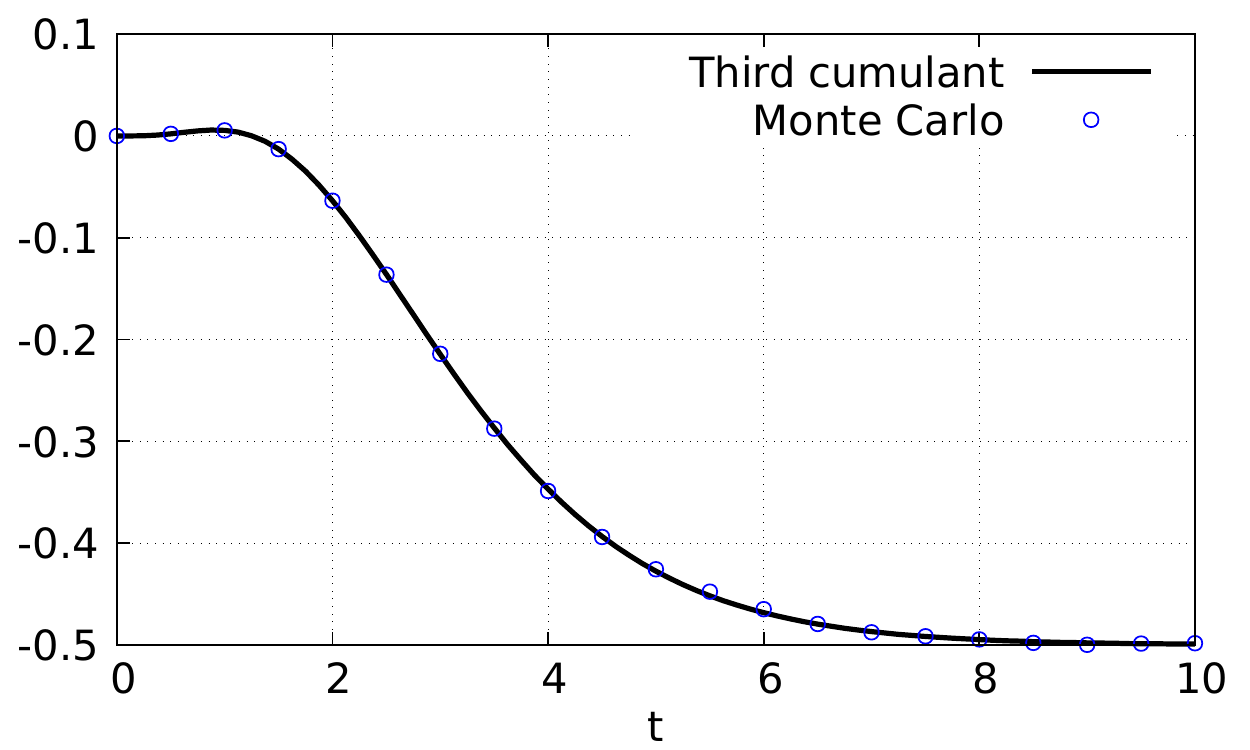}
    \caption{Third cumulant $\kappa^{(3)} (t)$.} 
 \end{subfigure}
 \begin{subfigure}[b]{0.49\textwidth}
    \includegraphics[width=1\linewidth, height=5cm]{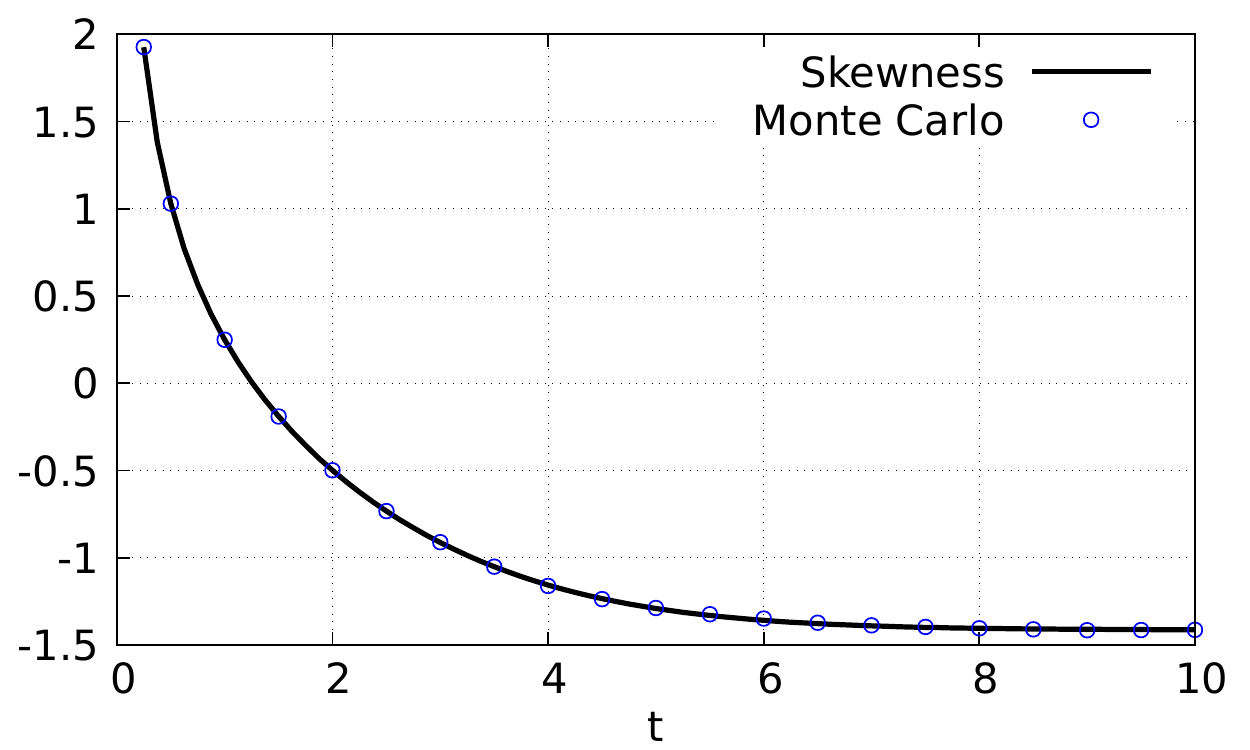}
 \caption{Skewness $\kappa^{(3)} (t)/\big( \kappa^{(2)} (t)\big)^{3/2}$.} 
 \end{subfigure}
 \caption{Third cumulant \eqref{c3} and skewness.} 
 \label{fig2} 
\end{figure}

In this section we consider the growth-collapse process 
$(X_t)_{t\in \real_+}$ of \cite{boxma}, defined as $X_t := t - Y_t$, i.e.
$$
X_t = t - \sum_{k=1}^{N_t} g(T_k,k,N_t) = t - \int_0^t g(s,N_s,N_t) dN_s,
\qquad t\in \real_+. 
$$
The moments of $X_t$
can be recovered from \eqref{eq4.1} 
and the binomial recursion 
\begin{equation}
  \label{fjk} 
\E [(X_t)^n]
=
\E [(t-Y_t)^n]
=
(-1)^n 
\left(
\E [ (Y_t)^n]
-
\sum_{k=0}^{n-1} {n \choose k} t^{n-k}
(-1)^k
\E [ (X_t)^k]
\right). 
\end{equation}
The next proposition extends Theorems~4 and 5 as well as Corollary~1
 of \cite{boxma} from mean and variance to moments
 of all orders, see also Corollary~4 in \cite{daw} for an expression
 using matrix exponentials. 
 It is also consistent with Theorem~3 of \cite{boxma}
 which states that the stationary distribution of the
  Markovian growth-collapse process $(X_t)_{t\in \real_+}$ is
  the gamma distribution $\Gamma (2,\lambda )$
  with shape parameter
  $2$ and scaling parameter $\lambda$, and
  cumulants $\kappa^{(n)}(\infty) = 2 (n-1)!/\lambda^n$, $n\geq 1$. 
 \begin{prop}
  \label{p5.1}
  The moments of the growth-collapse process
  $$
X_t = t - \sum_{k=1}^{N_t} f_k (T_k) (1-U_k) \prod_{l=k+1}^{N_t} U_l,
 \qquad t\in \real_+, 
$$
 with uniform cut-offs $(U_k)_{k\geq 1}$ on $[0,1]$, 
 are given by 
$$ 
  \E [ (X_t)^n ] = \frac{(n+1)!}{\lambda^n} \sum_{k=0}^n (-1)^k (k+1)^{n-1}{n \choose k}
  \re^{-k\lambda t/(k+1)},
 \quad n\geq 0, \quad t\in \real_+. 
$$ 
 As a consequence, the asymptotic moments of 
  $(X_t)_{t\in \real_+}$ are given by 
  $$
  \lim_{t\to \infty} 
  \E [ (X_t)^n ] = \frac{(n+1)!}{\lambda^n}, 
  \qquad n\geq 1. 
  $$
\end{prop}
Before proving Proposition~\ref{p5.1},
we recover the first moments and cumulants of $X_t$
 from the expressions \eqref{m1x}-\eqref{m4x}
 and the identity \eqref{fjk}.
 We find 
 \begin{equation}
   \label{abc1}
 \E [ X_t ] = - \E[Y_t] + t =  
 \frac{\lambda }{2} \re^{ - \lambda t / 2}  \int_0^t s_1
 \re^{  \lambda s_1 / 2} ds_1 
 = \frac{2}{\lambda} ( - \re^{-\lambda t/2} + 1 ) 
\end{equation}
 and 
\begin{equation}
\E[(X_t)^2] 
= \E [ (Y_t)^2] -t^2 + 2 t \E [ X_t ] 
= \frac{3!}{\lambda^2} \left( 3 \re^{ - 2 \lambda t / 3 } - 4 \re^{-\lambda t/2} + 1 \right), 
\end{equation}
see Theorems~4 and 5 of \cite{boxma}.
Next, from \eqref{m3x} we have 
\begin{eqnarray}
  \nonumber
  \E [ (X_t)^3] 
  & = &
  - \E [ (Y_t)^3]
  + t^3 - 3t^2 \E [ X_t]
  + 3 t \E [ X_t^2]
 \\
 & = &
   \frac{4!}{\lambda^3}
   \big( - 16 \re^{-3 \lambda t/  4} + 27 \re^{-2\lambda t/3}
   - 12 \re^{-\lambda t/2} + 1 
   \big), 
   \end{eqnarray} 
and therefore the third cumulant of $X_t$ is given by
\eqref{inv} as 
  \begin{equation}
    \label{c3}
\kappa^{(3)} (t) = 2 \frac{2!}{\lambda^3} 
\big( - 27 \re^{-7 \lambda t/6}
+ 96 \re^{-3 \lambda t/4}
-135 \re^{-2 \lambda t/3}
    + 4 \re^{-3 \lambda t/2}
    + 24 \re^{-\lambda t}
    + 39 \re^{-\lambda t/2}
    -1
    \big), \quad
\end{equation} 
  see Figure~\ref{fig2}.
   
\begin{figure}[H]
  \centering
 \begin{subfigure}[b]{0.49\textwidth}
    \includegraphics[width=1\linewidth, height=5cm]{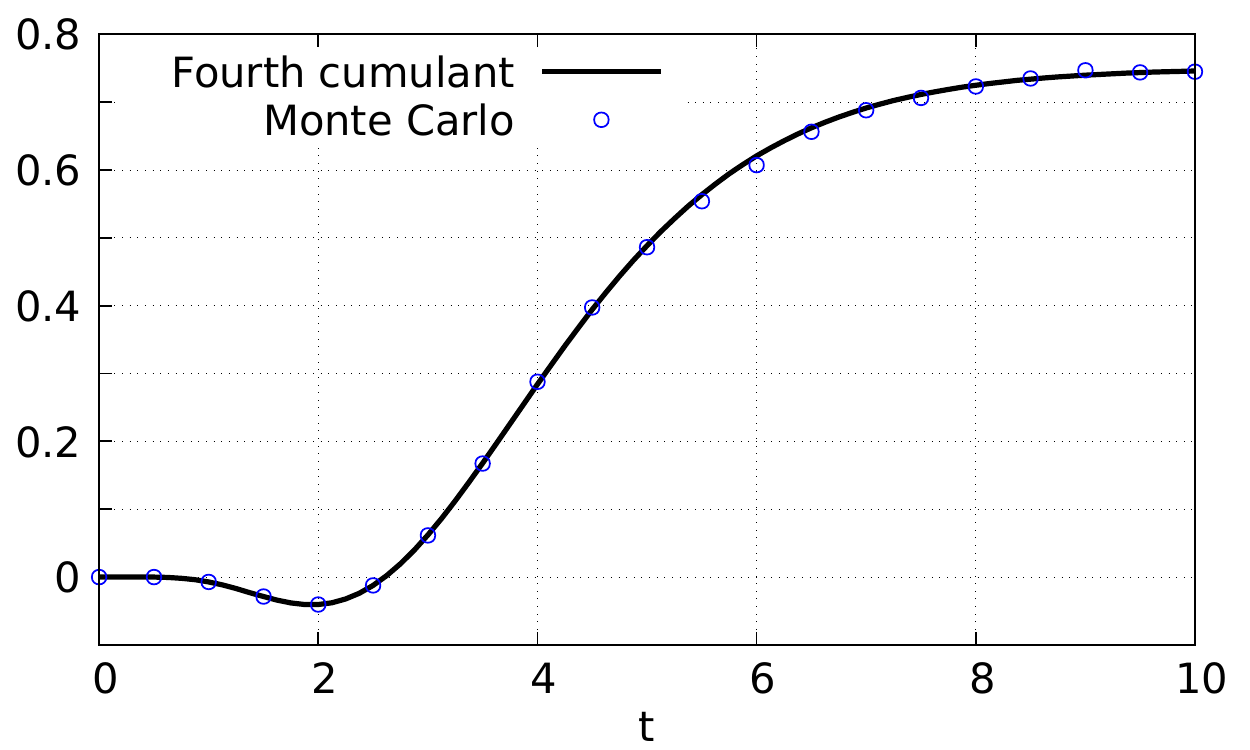}
    \caption{Fourth cumulant $\kappa^{(4)} (t)$.} 
 \end{subfigure}
 \begin{subfigure}[b]{0.49\textwidth}
    \includegraphics[width=1\linewidth, height=5cm]{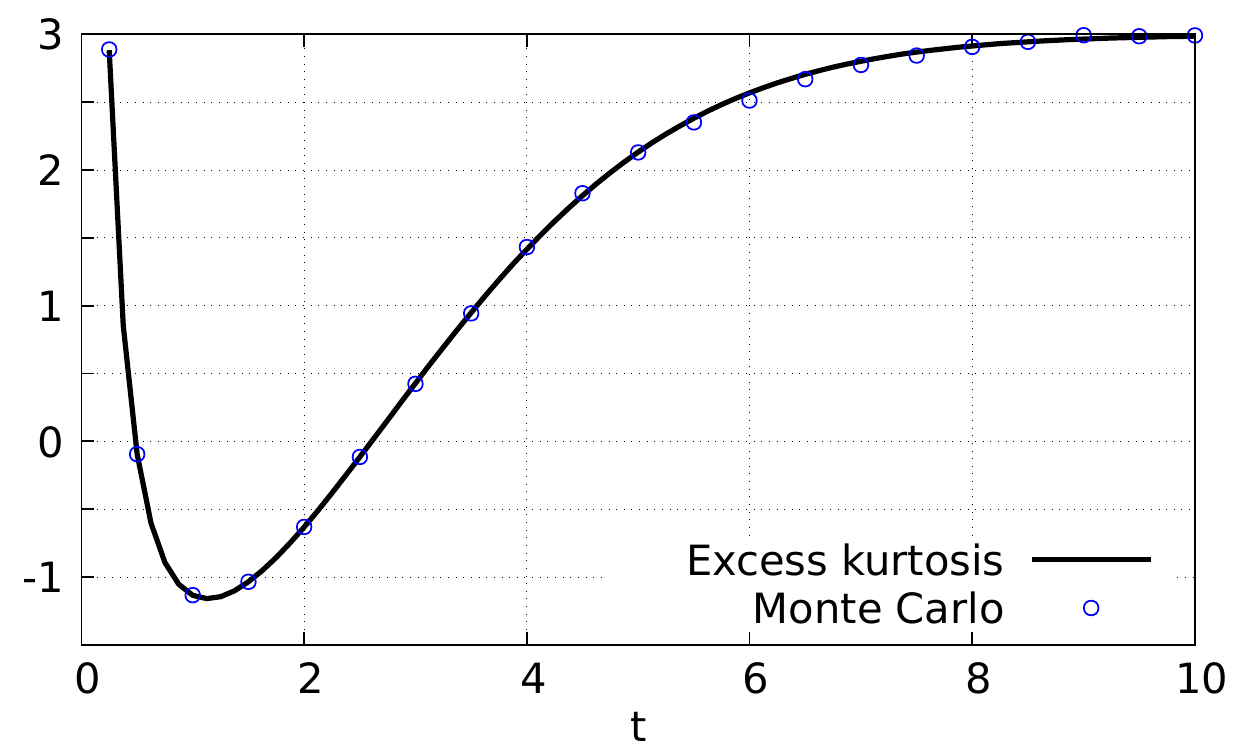} 
    \caption{Excess kurtosis $\kappa^{(4)} (t) / \big( \kappa^{(2)} (t) \big)^2$.} 
 \end{subfigure}
  \caption{Fourth cumulant \eqref{c4} and excess kurtosis.} 
\label{fig3} 
\end{figure}

\noindent
 From \eqref{m4x} and \eqref{abc1}-\eqref{c3} we have 
 \begin{eqnarray}
   \nonumber 
  \E [ (X_t)^4] 
  & = & \E [ (Y_t)^4]
  - t^4 + 4 t^3 \E [ X_t] - 6 t^2 \E [ (X_t)^2 ]
  + 4 t \E [ X_t^3]
  \\
  \label{abc4}
  & = & 
  \frac{5!}{\lambda^4} \big(
  125\re^{-4 \lambda t/5} 
 - 256 \re^{-3\lambda t/4}
 + 162 \re^{-2 \lambda t/3} 
 -32 \re^{- \lambda t/2 }
 + 1  \big), 
 \end{eqnarray} 
 and therefore the fourth cumulant of $X_t$ is given
 by \eqref{inv} as 
 \begin{eqnarray}
   \label{c4}
       \kappa^{(4)} (t) & =  & 
2 \frac{3!}{\lambda^4}
\big( 
   -8\re^{ -2\lambda t}
   + 504 \re^{- 7\lambda t/6}
   + 1250 \re^{ -4\lambda t/5}
   - 256 \re^{- 5\lambda t/4}
   - 2304 \re^{ -3\lambda t/4}
   \\
   \nonumber
   & &
   + 72 \re^{-5\lambda t/3}
   - 81 \re^{- 4\lambda t/3}
   + 1206 \re^{ -2\lambda t/3}
   - 64 \re^{-3\lambda t/2}
   - 168 \re^{ -\lambda t}
   - 152 \re^{ -\lambda t/2}
   + 1\big),
\end{eqnarray} 
 see Figure~\ref{fig3}. 
   
\begin{Proofy}{\em \ref{p5.1}}.
  Using the infinitesimal generator of the Markov process $(X_t)_{t\in \real_+}$ 
  it can be shown, see \S~3 of \cite{boxma} and \S~3.5 of \cite{daw}, that the moments 
  $\E [ (X_t)^n ]$ satisfy the differential equation
  $$
  \frac{d}{dt} \E [ (X_t)^n ] = n \E [ (X_t)^{n-1} ] - \frac{\lambda n}{n+1} \E [ (X_t)^n ]. 
$$ 
  Based on the intuition gained from 
 \eqref{abc1}-\eqref{abc4}, we search for a solution of the form
$$ 
 \E [ (X_t)^n ] = \frac{(n+1)!}{\lambda^n} \sum_{k=0}^n a_{k,n} \re^{-k\lambda t/(k+1)}, 
$$
 which, by identification of terms, yields the recurrence relation   
$$
 a_{k,n} = \frac{n(k+1)}{n-k} a_{k,n-1}, \qquad
 0 \leq k < n,
$$
 hence
$$
 a_{k,n} = (k+1)^{n-k} {n \choose k} a_{k,k}, \qquad
 0 \leq k < n.  
$$
 In addition, the initial condition $t=0$ requires 
$$ 
 \sum_{k=0}^n a_{k,n} = 0, 
$$
 hence
$$ 
 \sum_{k=0}^n (k+1)^{n-k} {n \choose k} a_{k,k} = 0, 
$$
 which is solved by taking $a_{k,k} = (-1)^k (k+1)^{k-1}$, 
 due to the combinatorial relation 
$$
 S(n,n+1) 
 = 
 \sum_{k=0}^n (-1)^{n-k} (k+1)^n
  {n+1 \choose k+1} 
 = 
 (n+1) \sum_{k=0}^n (k+1)^{n-1} {n \choose k} (-1)^k = 0, 
$$
 which follows from the vanishing of the Stirling numbers
 of the second kind $S(n,n+1)$, see e.g.
 page 824 of \cite{Hand1972}. 
\end{Proofy}
\section{Embedded growth-collapse chain}
\label{s6}
\noindent
In this section we show that Proposition~\ref{fldsf}
can also be used to compute the moments of all orders 
of the embedded chain
\begin{equation}
  \label{fjdksl} 
Y(m) = Y_{T_m}
=
\sum_{k=1}^m g(T_k,k,m) 
=
\int_0^\infty g(s,N_s,m) {\bf 1}_{[0,T_m]}(s) dN_s,
\qquad m \geq 1.
\end{equation} 

\begin{prop}
  \label{p6.1}
   Let $(Y(m))_{m\geq 1}$ be of the form
  \eqref{fjdksl}. 
  For any $n, m \geq 1$ we have 
  \begin{align*}
    & 
    \E \big[ ( Y (m) )^n \big] 
    =
  \\ 
\nonumber 
 &
 n! \sum_{k=1}^n
  \lambda^k
 \sum_{0=q_0 < q_1< \cdots < q_{k-1} < q_k = n }
 \E \left[ 
 \int_0^\infty \int_0^{s_k} \cdots \int_0^{s_2} 
 {\bf 1}_{\{N_{{s_k}^-} \leq m -k\}}
 \prod_{l=1}^k 
 \frac{\big( g(s_l,l+N_{s_l},m)
 \big)^{q_l-q_{l-1}}}{(q_l-q_{l-1})!} 
 d s_1 \cdots d s_k 
 \right]
.
\end{align*} 
\end{prop}
\begin{Proof}
 By Proposition~\ref{fldsf} and the identity
  $\{ s \leq T_m \} = \{ N_{s^-}< m\}$, $s >0$, for all $n\geq 1$ we have
\begin{align*} 
 & 
 \E \left[ 
 \left(
\sum_{k=1}^m g(T_k,k,m) 
 \right)^n 
 \right] 
 =
 \E \left[ 
 \left(
 \int_0^\infty g(s,N_s,m) {\bf 1}_{[0,T_m]}(s) dN_s
 \right)^n 
 \right] 
\\ 
\nonumber 
 & =  
 \sum_{k=1}^n
 k!
 \lambda^k
 \sum_{ \pi_1 \cup \cdots \cup \pi_k = \{ 1 , \ldots , n\} } 
 \E \left[ 
 \int_0^\infty \int_0^{s_k} \cdots \int_0^{s_2} 
 \epsilon^+_{s_1} \cdots \epsilon^+_{s_k} 
 \prod_{l=1}^k 
 \big( g(s_l,N_{s_l},m)
 {\bf 1}_{\{N_{s_l^-}< m\}}
 \big)^{|\pi_l|} 
 d s_1 \cdots d s_k 
 \right]
\\ 
\nonumber 
 & =  
 \sum_{k=1}^n
 k!
 \lambda^k
 \sum_{ \pi_1 \cup \cdots \cup \pi_k = \{ 1 , \ldots , n\} } 
 \E \left[ 
 \int_0^\infty \int_0^{s_k} \cdots \int_0^{s_2} 
 \prod_{l=1}^k 
 \big( g(s_l,l+N_{s_l},m)
 {\bf 1}_{\{l-1+N_{s_l^-}< m\}}
 \big)^{|\pi_l|} 
 d s_1 \cdots d s_k 
 \right]
\\ 
\nonumber 
 & =  
 \sum_{k=1}^n
 k!
 \lambda^k
 \sum_{ \pi_1 \cup \cdots \cup \pi_k = \{ 1 , \ldots , n\} } 
 \E \left[ 
 \int_0^\infty \int_0^{s_k} \cdots \int_0^{s_2} 
 {\bf 1}_{\{N_{{s_k}^-} \leq m -k\}}
 \prod_{l=1}^k 
 \big( g(s_l,l+N_{s_l},m)
 \big)^{|\pi_l|} 
 d s_1 \cdots d s_k 
 \right]
.
\end{align*} 
\end{Proof}
 Next, we specialize Proposition~\ref{p6.1} to the case where 
 $g(s,k,m)$ takes the
form
\begin{equation}
  \label{fjl} 
g(s,k,m) = f_k(s) (1-Z_k) \prod_{l=k+1}^m Z_l,
\end{equation} 
where $(Z_k)_{k \geq 1}$ is an i.i.d. random sequence independent
of the Poisson process $(N_t)_{t\in \real_+}$, with moment sequence
$m_n = \E [ Z^n]$, $n\geq 0$. 
  
\begin{corollary} 
  \label{p6.2}
  Let $(Y(m))_{m\geq 1}$ be defined as in \eqref{fjl}. 
 For any $n, m \geq 1$ we have 
\begin{align} 
  \nonumber
 & \E \big[ ( Y (m) )^n \big] 
 = 
 n! \sum_{k=1}^n
 \sum_{i=0}^{m-k} \frac{\lambda^{k+i}}{i!}
 m_n^{m - i-k}
\sum_{0=q_0 < q_1< \cdots < q_{k-1} < q_k = n }
 \\ 
\nonumber 
 & 
\prod_{l=1}^k \frac{C_{q_{l-1},q_l-q_{l-1}}}{(q_l-q_{l-1})!}
\int_0^\infty \re^{-\lambda s_k} \int_0^{s_k} \cdots \int_0^{s_2} 
 \left(
 s_1 + \sum_{l=1}^{k-1}
 m_{p_1+\cdots + p_l} ( s_{l+1} - s_l ) 
 \right)^i
 \prod_{l=1}^k
 f_{l,m}^{q_l-q_{l-1}}(s_l) d s_1 \cdots d s_k.
 \\
 \label{fjkldsf3}
\end{align} 
\end{corollary}
\begin{Proof}
  Using the fact that the jumps of $(N_s)_{s\in [0,s_k]}$
  are uniformly distributed on $[0,s_k]$ given that
  $N_{s_k}=i$, by Proposition~\ref{p6.1} we have
\begin{align*} 
 & 
 \E \left[ 
 \left(
\sum_{k=1}^m g(T_k,k,m) 
 \right)^n 
 \right] 
 = \sum_{k=1}^n k! \lambda^k \sum_{ \pi_1 \cup \cdots \cup \pi_k = \{ 1 , \ldots , n\} }
 \\ \nonumber & \quad \E \left[ \int_0^\infty \int_0^{s_k} \cdots \int_0^{s_2} {\bf 1}_{\{N_{{s_k}^-} \leq m -k \}} \prod_{l=1}^k \left( f_{l,m}(s_l) (1-Z_{l+N_{s_l}} ) \prod_{j=l+1+N_{s_l}}^m Z_j \right)^{|\pi_l|} d s_1 \cdots d s_k \right]
\\ 
\nonumber 
 & =  
 \sum_{k=1}^n
 k!
 \lambda^k
 \sum_{ \pi_1 \cup \cdots \cup \pi_k = \{ 1 , \ldots , n\} } 
 \\
 \nonumber
 &
 \quad \E \left[ 
 \int_0^\infty \int_0^{s_k} \cdots \int_0^{s_2} 
 {\bf 1}_{\{N_{{s_k}^-} \leq m -k \}}
 \prod_{l=1}^k
 \left(
 f_{l,m}^{p_l}(s_l) C_{p_1+\cdots + p_{l-1},p_l}
 ( m_{p_1+\cdots + p_l} )^{N_{s_{l+1}} - N_{s_l}}
 \right)
 d s_1 \cdots d s_k 
 \right]
\\ 
\nonumber 
 & =  
 n!\sum_{k=1}^n
 \frac{ \lambda^k}{p_1!\cdots p_k!}
 \sum_{i=0}^{m-k} \P ( N_{{s_k}^-} = i ) 
 \sum_{p_1 + \cdots + p_k = n \atop p_1,\ldots, p_k \geq 1} 
 \\
 & \quad 
 \int_0^\infty \int_0^{s_k} \cdots \int_0^{s_2} 
 \prod_{l=1}^k
 \big( f_{l,m}^{p_l}(s_l) C_{p_1+\cdots + p_{l-1},p_l} \big) 
 \left(
 \frac{s_1 + \sum_{l=1}^{k-1}
 m_{p_1+\cdots + p_l} ( s_{l+1} - s_l ) 
 }{s_k} \right)^i
 m_n^{m - i-k}
 d s_1 \cdots d s_k 
\\ 
\nonumber 
 & =  
 n! \sum_{k=1}^n
 \lambda^{k+i}
 \sum_{p_1 + \cdots + p_k = n \atop p_1,\ldots, p_k \geq 1} 
 \sum_{i=0}^{m-k} \frac{m_n^{m - i-k}}{i!}
 \prod_{l=1}^k
 \frac{C_{p_1+\cdots + p_{l-1},p_l}}{p_l!}
\\
 & \quad  
 \int_0^\infty \re^{-\lambda s_k} \int_0^{s_k} \cdots \int_0^{s_2} 
 \prod_{l=1}^k
 f_{l,m}^{p_l}(s_l)
 \left(
 s_1 + \sum_{l=1}^{k-1}
 m_{p_1+\cdots + p_l} ( s_{l+1} - s_l ) 
 \right)^i
 d s_1 \cdots d s_k
.
\end{align*} 
\end{Proof}
\subsubsection*{Uniform cut-offs} 
The next result specializes Corollary~\ref{p6.2}
to the case of embedded growth processes with uniform
cut-offs.
\begin{corollary}
  \label{p6.3}
  Let $(Y(m))_{m\geq 1}$ be
    defined as in \eqref{fjl}, 
    with $(Z_k)_{k \geq 1}$ an i.i.d. uniform random sequence.
    For any $n, m \geq 1$ we have 
\begin{eqnarray} 
  \nonumber
    \E \big[ ( Y (m) )^n \big] 
 & = & 
 \sum_{k=1}^n
 \sum_{i=0}^{m-k} \frac{\lambda^{k+i}}{i!(n+1)^{m - i-k}}
 \sum_{0=q_0 < q_1< \cdots < q_{k-1} < q_k = n }
 \prod_{l=1}^k \frac{1}{1+q_l}
  \\ 
\nonumber 
 & &
 \quad \int_0^\infty \re^{-\lambda s_k} \int_0^{s_k} \cdots \int_0^{s_2} 
 \left(
 s_1 + \sum_{l=1}^{k-1}
 \frac{s_{l+1} - s_l}{1+q_l} 
 \right)^i
 \prod_{l=1}^k
 f_{l,m}^{q_l-q_{l-1}}(s_l)
 d s_1 \cdots d s_k.
 \\
 \label{fjkldsf3-2}
\end{eqnarray} 
\end{corollary}
\begin{Proof}
 By Corollary~\ref{p6.2} we have
\begin{eqnarray*} 
\lefteqn{ 
 \E \left[ 
 \left(
\sum_{k=1}^m g(T_k,k,m) 
 \right)^n 
 \right] 
 =  
 n! \sum_{k=1}^n
 \lambda^{k+i}
 \sum_{p_1 + \cdots + p_k = n \atop p_1,\ldots, p_k \geq 1} 
 \sum_{i=0}^{m-k} \frac{m_n^{m - i-k}}{i!}
 \prod_{l=1}^k
 \frac{C_{p_1+\cdots + p_{l-1},p_l}}{p_l!}
}
\\
 & & 
 \int_0^\infty \re^{-\lambda s_k} \int_0^{s_k} \cdots \int_0^{s_2} 
 \prod_{l=1}^k
 f_{l,m}^{p_l}(s_l)
 \left(
 s_1 + \sum_{l=1}^{k-1}
 m_{p_1+\cdots + p_l} ( s_{l+1} - s_l ) 
 \right)^i
 d s_1 \cdots d s_k
 \\ 
\nonumber 
 & = &   
 n! \sum_{k=1}^n
 \lambda^{k+i}
 \sum_{p_1 + \cdots + p_k = n \atop p_1,\ldots, p_k \geq 1} 
 \sum_{i=0}^{m-k} \frac{m_n^{m - i-k}}{i!}
 \prod_{l=1}^k
 \frac{(p_1+\cdots + p_{l-1})!}{(p_1+\cdots + p_l+1)!}
\\
 & & 
 \int_0^\infty \re^{-\lambda s_k} \int_0^{s_k} \cdots \int_0^{s_2} 
 \prod_{l=1}^k
 f_{l,m}^{p_l}(s_l)
 \left(
 s_1 + \sum_{l=1}^{k-1}
 m_{p_1+\cdots + p_l} ( s_{l+1} - s_l ) 
 \right)^i
 d s_1 \cdots d s_k. 
\end{eqnarray*} 
\end{Proof}
The first moments of $Y(m)$ can be computed in closed form when $f_{l,m}(s)=s$,
which corresponds to 
$$ 
Y(m) = Y_{T_m}
=
\sum_{k=1}^m T_k (1-U_k) \prod_{l=k+1}^m U_l,
\qquad m \geq 1, 
$$
 where $(U_k)_{k\geq 1}$ is a uniform random sequence on $[0,1]$. 
\subsubsection*{First moment} 
 For $n=1$, Corollary~\ref{p6.3} yields 
$$
\E \big[ Y(m) \big]
=  \frac{1}{\lambda}
\sum_{i=0}^{m-1} 
    \frac{1}{i!2^{m - i-1}}
    \int_0^\infty 
 \re^{- s_1} s_1^i \frac{s_1}{2}
 d s_1 
=
 \frac{2^{-m}+m-1}{\lambda}
, 
$$ 
 which is consistent with Theorem~7 in \cite{boxma}.
\subsubsection*{Second moment} 
 For $n=2$, we find 
\begin{align*} 
 & 
  \E \big[ ( Y(m) )^2 \big]
\\
 & = 
\frac{1}{\lambda^2}
\sum_{i=0}^{m-1} 
 \frac{1}{i!3^{m-1-i}}
 \int_0^\infty 
 \re^{- s_1} s_1^i \frac{s_1^2}{3}
 d s_1 
 +  
 \frac{1}{2\lambda^2}
 \sum_{i=0}^{m-2} 
 \frac{1}{i! 3^{m-2-i}}
 \int_0^\infty
 \re^{- s_2} 
 \int_0^{s_2} 
 \left(\frac{s_1+s_2}{2} \right)^i 
 s_1
 \frac{s_2}{3}
 d s_1 d s_2
 \\
 & =  
 \frac{1}{\lambda^2}
 \left(
 \frac{2}{3^m} 
 +
 \frac{m-1}{2^{m-1}} -m + m^2
 \right). 
\end{align*} 
\subsubsection*{Third moment} 
\begin{figure}[H]
  \centering
 \begin{subfigure}[b]{0.49\textwidth}
    \includegraphics[width=1\linewidth, height=5cm]{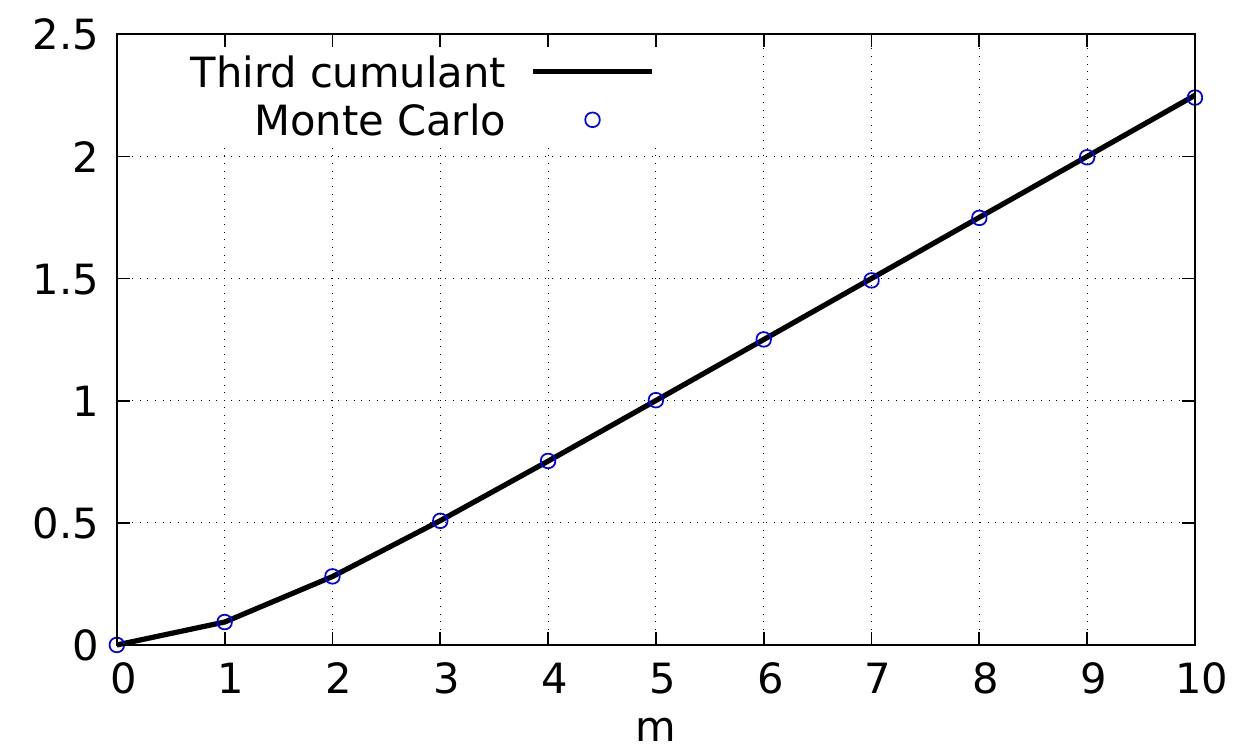}
    \caption{Third cumulant of $Y(m)$.} 
 \end{subfigure}
 \begin{subfigure}[b]{0.49\textwidth}
    \includegraphics[width=1\linewidth, height=5cm]{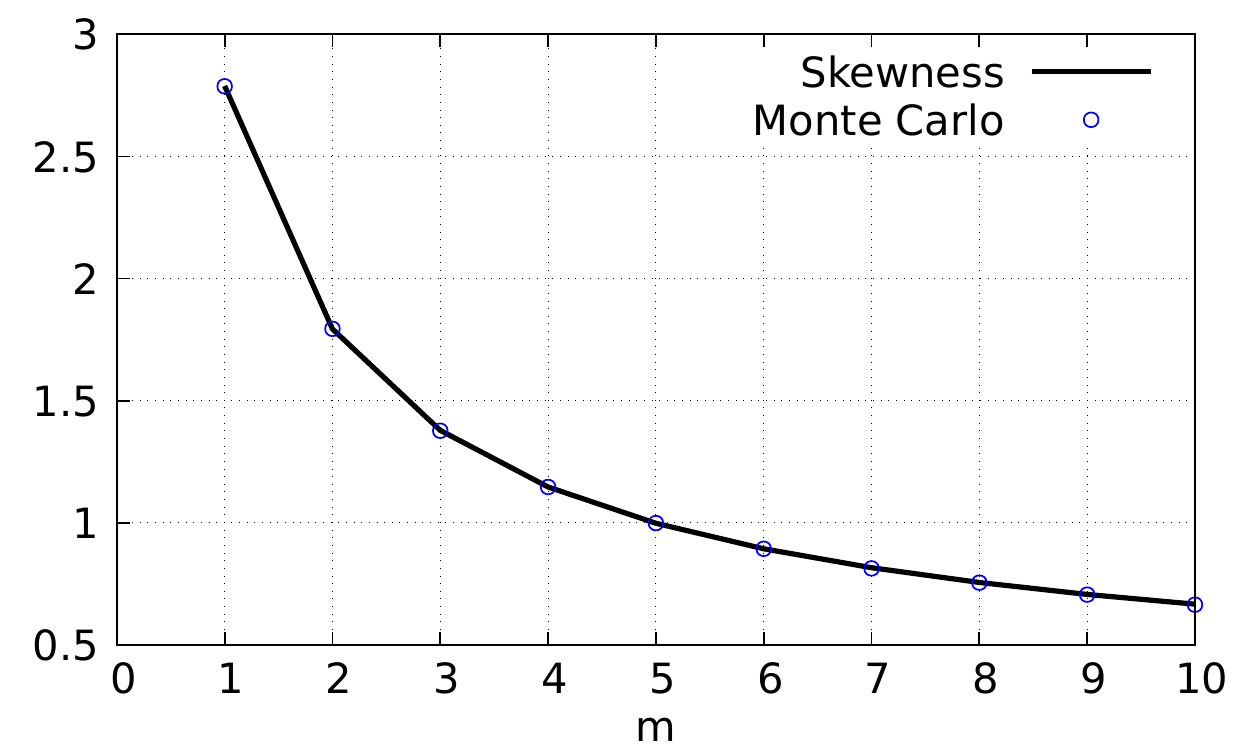} 
 \caption{Skewness $\kappa^{(3)} (m)/\big( \kappa^{(2)} (m)\big)^{3/2}$.} 
 \end{subfigure}
 \caption{Third cumulant and skewness of $Y(m)$.} 
 \label{fig2.1} 
\end{figure}

 For $n=3$, we have 
 \begin{eqnarray}
   \nonumber 
  \E \big[ Y(m)^3 \big] 
 & = & 
  \frac{1}{\lambda^3}
  \sum_{i=0}^{m-1} \frac{1}{i!4^{m-1-i}}  
 \int_0^\infty \re^{- s_1}
 \frac{s_1^{i+3}}{4}
 ds_1 
 \\
 \nonumber 
 & &
 + 
 \frac{1}{\lambda^3}
 \sum_{i=0}^{m-2} \frac{1}{i!4^{m-2-i}}
 \int_0^\infty \re^{- s_2} \int_0^{s_2} 
 \left(
 \left(
 \frac{s_1+s_2}{2} 
 \right)^i
 \frac{s_1}{2}
 \frac{s_2^2}{4} 
 +
 \left(
 \frac{2s_1 + s_2}{3} 
 \right)^i
 \frac{s_1^2}{3}
 \frac{s_2}{4}
 \right) d s_1 d s_2
 \\
\label{fkjlds1} 
 & &
 +
 \frac{1}{\lambda^3}
 \sum_{i=0}^{m-3} \frac{1}{i!4^{m-3-i}}
 \int_0^\infty \re^{- s_3} \int_0^{s_3} \int_0^{s_2} 
 \left(
 \frac{3s_1 + s_2 + 2 s_3}{6} 
 \right)^i
 \frac{s_1}{2}
 \frac{s_2}{3}
 \frac{s_3}{4}
 d s_1 d s_2 ds_3 . 
\end{eqnarray} 
Although the last partial summation \eqref{fkjlds1} does not have a polynomial
expression, it can easily be estimated numerically, see Figure~\ref{fig2.1}.

\subsubsection*{Fourth moment} 
\begin{figure}[H]
  \centering
 \begin{subfigure}[b]{0.49\textwidth}
    \includegraphics[width=1\linewidth, height=5cm]{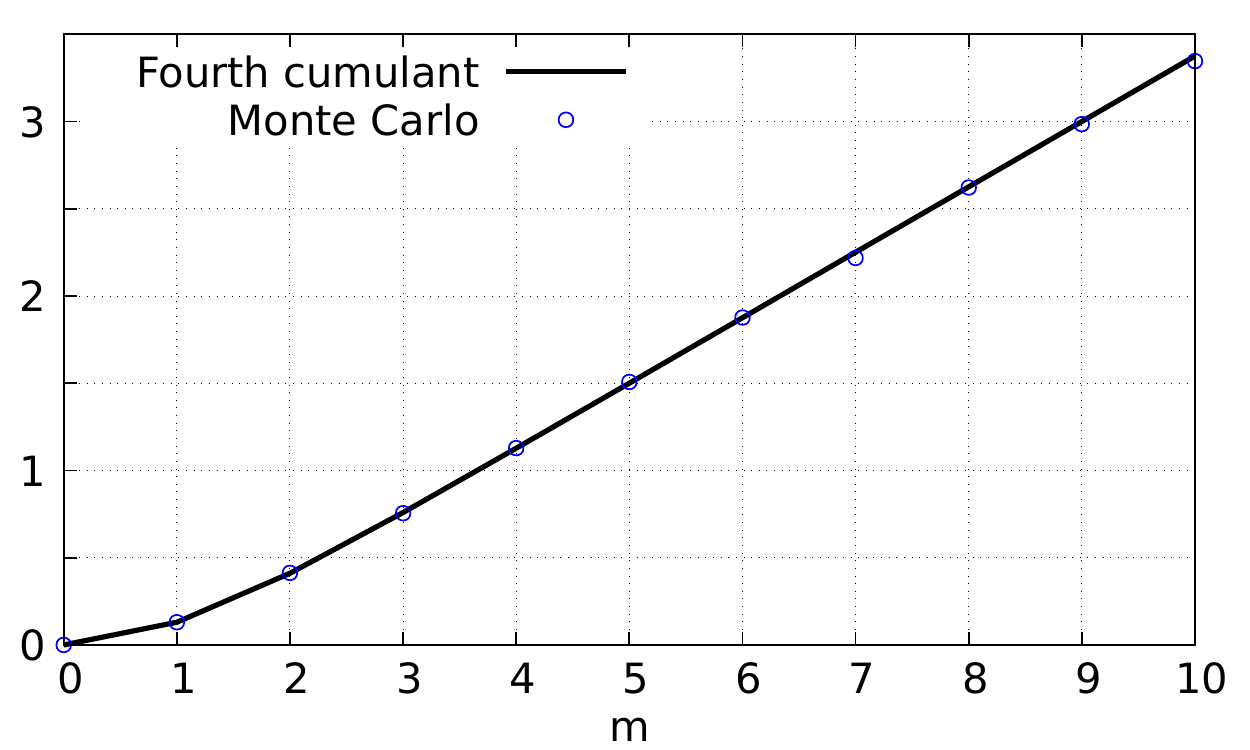}
    \caption{Fourth cumulant of $Y(m)$.} 
 \end{subfigure}
 \begin{subfigure}[b]{0.49\textwidth}
    \includegraphics[width=1\linewidth, height=5cm]{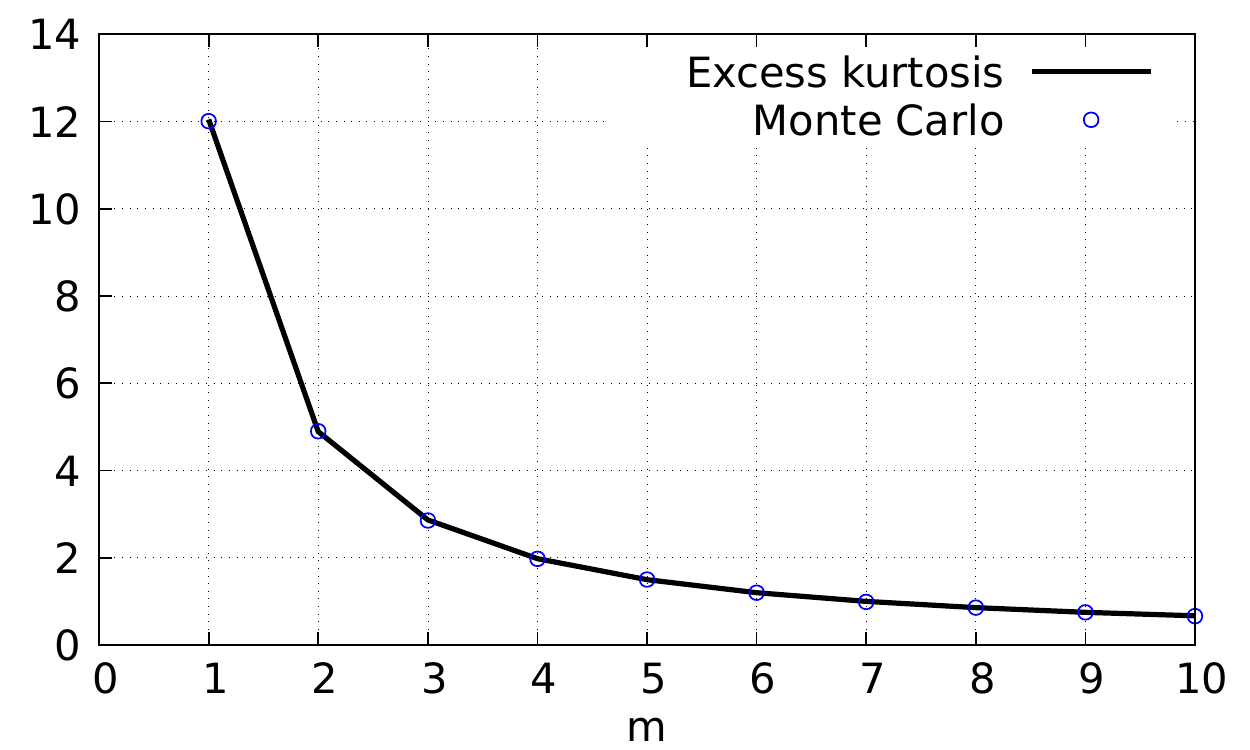} 
    \caption{Excess kurtosis $\kappa^{(4)} (t) / \big( \kappa^{(2)} (m) \big)^2$.} 
 \end{subfigure}
 \caption{Fourth cumulant and excess kurtosis of $Y(m)$.} 
\label{fig3.1} 
\end{figure}

  \noindent 
 For $n=4$, we find 
\begin{align*} 
 & 
    \E \big[ Y(m)^4 \big] 
 = 
 \frac{1}{\lambda^4}
 \sum_{i=0}^{m-1} \frac{1}{i!5^{m-1-i}}
 \int_0^\infty \re^{- s_1} \frac{s_1^{i+4}}{5}
 d s_1 
    \\
  &  +  \frac{1}{\lambda^4}
  \sum_{i=0}^{m-2} \frac{1}{i!5^{m-2-i}}
  \int_0^\infty \re^{- s_2} \int_0^{s_2}
  \left(
 \left(
 \frac{s_1+s_2}{2} 
 \right)^i
  \frac{s_1}{2}
 \frac{s_2^3}{5}
 + \left(
 \frac{3s_1 + s_2}{4} 
 \right)^i
 \frac{s_1^3}{4}
 \frac{s_2}{5}
 +
 \left(
 \frac{2s_1 + s_2}{3} 
 \right)^i
 \frac{s_1^2}{3}
 \frac{s_2^2}{5}
 \right)
 d s_1 d s_2
  \\
  &  + 
  \frac{1}{\lambda^4}
  \sum_{i=0}^{m-3} \frac{1 }{i!5^{m-3-i}}
 \int_0^\infty \re^{- s_3} \int_0^{s_3} \int_0^{s_2}
 \\
  & 
 \left(
 \left(
 \frac{8s_1 + s_2 + 3 s_3}{12} 
 \right)^i
 \frac{s_1^2}{3}
 \frac{s_2}{4}
 \frac{s_3}{5}
 +
 \left(
 \frac{4s_1 + 2s_2 + 2 s_3}{8} 
 \right)^i
 \frac{s_1}{2}
 \frac{s_2^2}{4}
 \frac{s_3}{5}
 +
 \left(
 \frac{3s_1 + s_2 + 2 s_3}{6} 
 \right)^i
 \frac{s_1}{2}
 \frac{s_2}{3}
 \frac{s_3^2}{5}
 \right)
 d s_1 ds_2 d s_3
 \\
  &  + 
 \frac{1}{\lambda^4}
 \sum_{i=0}^{m-4} \frac{1}{i!5^{m-4-i}}
 \int_0^\infty \re^{- s_4} \int_0^{s_4} \int_0^{s_3} \int_0^{s_2} 
 \left(
 \frac{6s_1 + 2 s_2 + s_3 + 3 s_4}{12} 
 \right)^i
 \prod_{l=1}^4
 \frac{s_1}{2}
 \frac{s_2}{3}
 \frac{s_3}{4}
 \frac{s_4}{5}
 d s_1 ds_2 ds_3 s_4, 
\end{align*} 
see Figure~\ref{fig3.1}. 

\subsubsection*{Compensated embedded chain}
\noindent
 Finally, we consider the compensated embedded chain 
 \begin{eqnarray}
   \nonumber 
 X(m) & := & Y_{T_m} - f_m(T_m)
 \\
 \nonumber
 & = & 
 - f_m(T_m) + \sum_{k=1}^m g(T_k,k,m) 
 \\
 \label{fjk3}
 & = & 
 - f_m(T_m) + \sum_{k=1}^m f_k(s) (1-Z_k) \prod_{l=k+1}^n Z_l,
 \qquad
 m \geq 1, 
\end{eqnarray} 
 where $(Z_k)_{k \geq 1}$ is an i.i.d. random sequence independent
 with moment sequence $m_n = \E [ Z^n]$, $n\geq 0$. 
 This process can be obtained from 
 \eqref{fjdksl} by replacing 
 $g(s,m,m) = f_m(s) (1-Z_m)$
 with $g(s,m,m) = f_m(s) (1-Z_m) - f_m(s)= - f_m(s) Z_m$ in \eqref{fjdksl}.
 In Corollary~\ref{p6.2}, this amounts to
 modifying last term ofr $i=m-k$ or $N_{s_k}=m-k$
 in \eqref{fjkldsf3},
 by changing 
 the last term $C_{q_{k-1},q_k-q_{k-1}}$ of order $l=k$ in the product
 $\displaystyle \prod_{l=1}^k C_{q_{l-1},q_l-q_{l-1}}$
 into 
 $(-1)^{q_k-q_{k-1}}m_{q_k} =(-1)^{q_k-q_{k-1}} m_n$, 
 yielding the next corollary. 
 \begin{corollary}
   \label{p6.4}
 Let $(X(m))_{m\geq 1}$ be defined as in \eqref{fjk3}. 
 For any $n, m \geq 1$ we have 
\begin{align} 
  \nonumber
 & \E \big[ ( X (m) )^n \big] 
 = 
 n! \sum_{k=1}^n
 \sum_{i=0}^{m-k-1} \frac{\lambda^{k+i}}{i!}
 m_n^{m - i-k}
\sum_{0=q_0 < q_1< \cdots < q_{k-1} < q_k = n }
 \\ 
\nonumber 
 & 
\int_0^\infty \re^{-\lambda s_k} \int_0^{s_k} \cdots \int_0^{s_2} 
 \left(
 s_1 + \sum_{l=1}^{k-1}
 m_{q_l} ( s_{l+1} - s_l ) 
 \right)^i
 \prod_{l=1}^k \frac{C_{q_{l-1},q_l-q_{l-1}} f_{l,m}^{q_l-q_{l-1}}(s_l)}{(q_l-q_{l-1})!}
 d s_1 \cdots d s_k
 \\
 \nonumber
 & + n! \sum_{k=1}^{\min ( n , m) } 
\frac{\lambda^m}{(m-k)!}
\sum_{0=q_0 < q_1< \cdots < q_{k-1} < q_k = n }
 (-1)^{q_k-q_{k-1}} m_n 
\\ 
\nonumber 
 & 
 \int_0^\infty \re^{-\lambda s_k}
f_{k,m}^{n-q_{k-1}}(s_k)
\int_0^{s_k} \cdots \int_0^{s_2} 
 \left(
 s_1 + \sum_{l=1}^{k-1}
 m_{q_l} ( s_{l+1} - s_l ) 
 \right)^{m-k}
 \prod_{l=1}^{k-1} \frac{C_{q_{l-1},q_l-q_{l-1}} f_{l,m}^{q_l-q_{l-1}}(s_l)}{(q_l-q_{l-1})!}
 d s_1 \cdots d s_k. 
\end{align} 
\end{corollary} 
  Similarly, according to Corollary~\ref{p6.3}, computing the moments of
\begin{equation}
  \label{fjk3-2}
 X(m) = T_m - \sum_{k=1}^m T_k (1-U_k) \prod_{l=k+1}^m U_l, 
\end{equation} 
 with $(U_k)_{k \geq 1}$ an i.i.d. uniform random sequence, 
 means multiplying the product
 $\displaystyle \prod_{l=1}^k \frac{1}{1+q_l}$
 for $i=m-k$ in \eqref{fjkldsf3-2}
 by 
 $$
 \frac{
   (-1)^{q_k-q_{k-1}}m_{q_k} 
 }{C_{q_{k-1},q_k-q_{k-1}}}
 =
 \frac{
 (-1)^{q_k-q_{k-1}}
 (1+q_k)! }{ (1+q_k) q_{k-1}!(q_k-q_{k-1})!}
 =
 (-1)^{q_k-q_{k-1}}
 {n \choose q_{k-1}}, 
 $$
 as is done in the next corollary. 
\begin{corollary} 
   \label{p6.5}
 Let $(X(m))_{m\geq 1}$ be defined as in \eqref{fjk3-2}, 
 with $(U_k)_{k \geq 1}$ an i.i.d. uniform random sequence.
 For any $n, m \geq 1$ we have 
\begin{eqnarray} 
  \nonumber
     \E \big[ ( Y (m) )^n \big] 
 & = & 
 \sum_{k=1}^n
 \sum_{i=0}^{m-k-1} \frac{\lambda^{k+i}}{i!(n+1)^{m - i-k}}
 \sum_{0=q_0 < q_1< \cdots < q_{k-1} < q_k = n }
  \\ 
\nonumber 
 & &
 \quad \int_0^\infty \re^{-\lambda s_k} \int_0^{s_k} \cdots \int_0^{s_2} 
 \left(
 s_1 + \sum_{l=1}^{k-1}
 \frac{s_{l+1} - s_l}{1+q_l} 
 \right)^i
 \prod_{l=1}^k \frac{f_{l,m}^{q_l-q_{l-1}}(s_l)}{1+q_l}
 d s_1 \cdots d s_k
 \\
 \nonumber
 & &
 + \sum_{k=1}^{\min (n,m )}
 \frac{\lambda^{m}}{(m-k)!}
 \sum_{0=q_0 < q_1< \cdots < q_{k-1} < q_k = n }
   (-1)^{q_k-q_{k-1}}
 {n \choose q_{k-1}}
 \\ 
\nonumber 
 & &
 \quad \int_0^\infty \re^{-\lambda s_k} \int_0^{s_k} \cdots \int_0^{s_2} 
 \left(
 s_1 + \sum_{l=1}^{k-1}
 \frac{s_{l+1} - s_l}{1+q_l} 
 \right)^{m-k}
\prod_{l=1}^k \frac{f_{l,m}^{q_l-q_{l-1}}(s_l)}{1+q_l}
   d s_1 \cdots d s_k.
\end{eqnarray} 
\end{corollary}
 When $f_{l,m}(s)=s$, Corollary~\ref{p6.5} shows that the second moment reads 
 $$
 \E \big[ ( X(m) )^2 \big]
 =
 \frac{1}{\lambda^2}
 \left(
 2
 - 4 \left(\frac{1}{2}\right)^m
 + 2 \left(\frac{1}{3}\right)^m
 \right) 
 $$
 which recovers Theorem~7 in \cite{boxma}.

 \begin{figure}[H]
  \centering
 \begin{subfigure}[b]{0.49\textwidth}
    \includegraphics[width=1\linewidth, height=5cm]{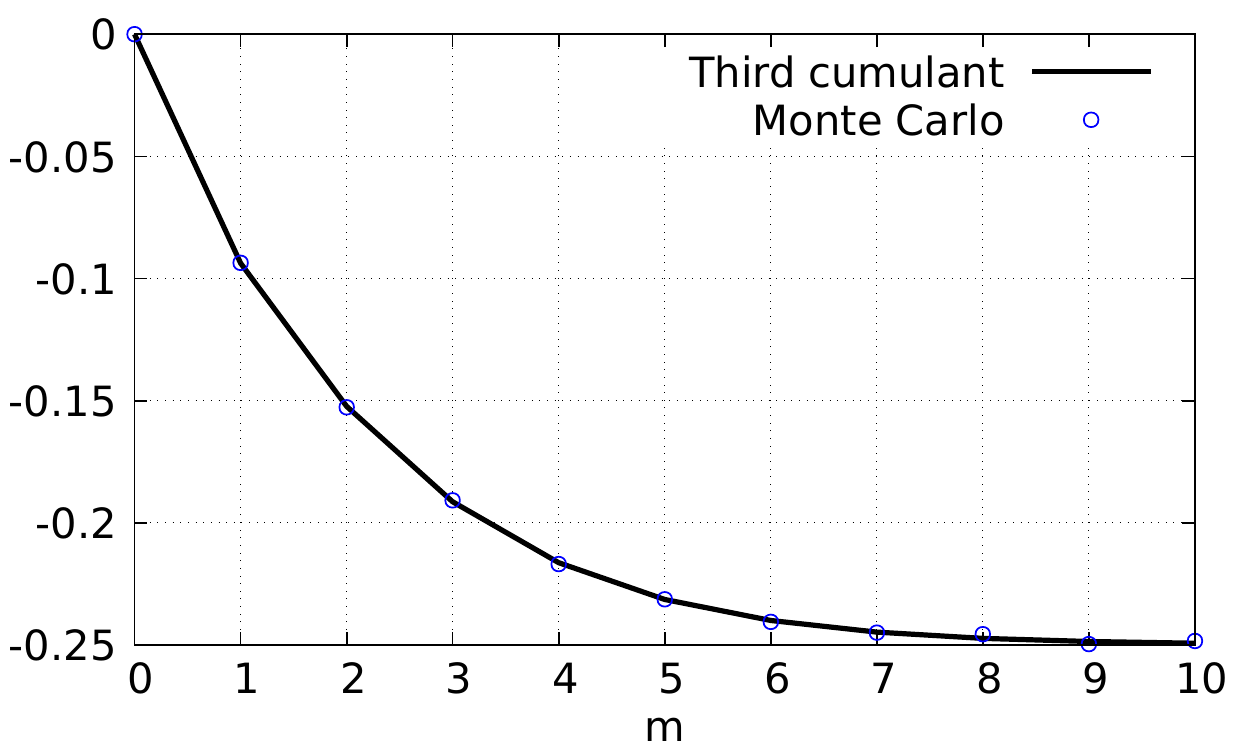}
    \caption{Third cumulant of $X(m)$.} 
 \end{subfigure}
 \begin{subfigure}[b]{0.49\textwidth}
    \includegraphics[width=1\linewidth, height=5cm]{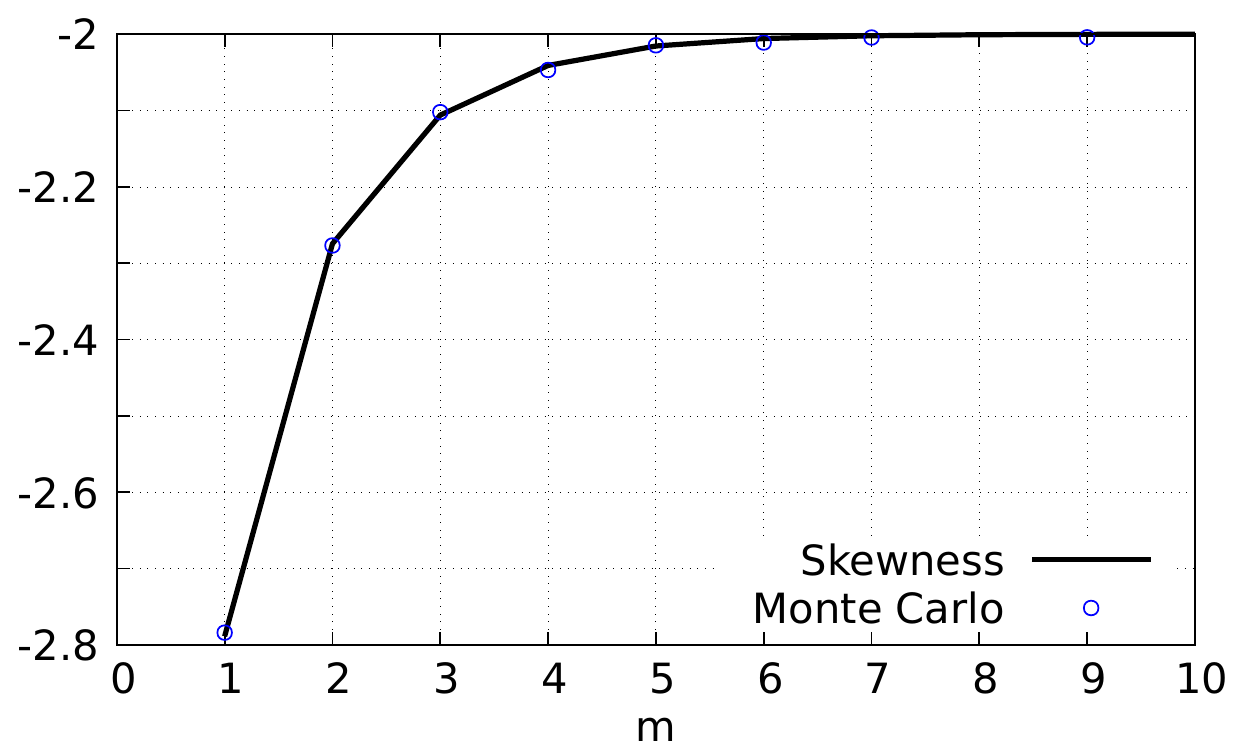} 
 \caption{Skewness $\kappa^{(3)} (m)/\big( \kappa^{(2)} (m)\big)^{3/2}$.} 
 \end{subfigure}
 \caption{Third cumulant and skewness of $X(m)$.} 
 \label{fig6} 
 \end{figure}
\begin{figure}[H]
  \centering
 \begin{subfigure}[b]{0.49\textwidth}
    \includegraphics[width=1\linewidth, height=5cm]{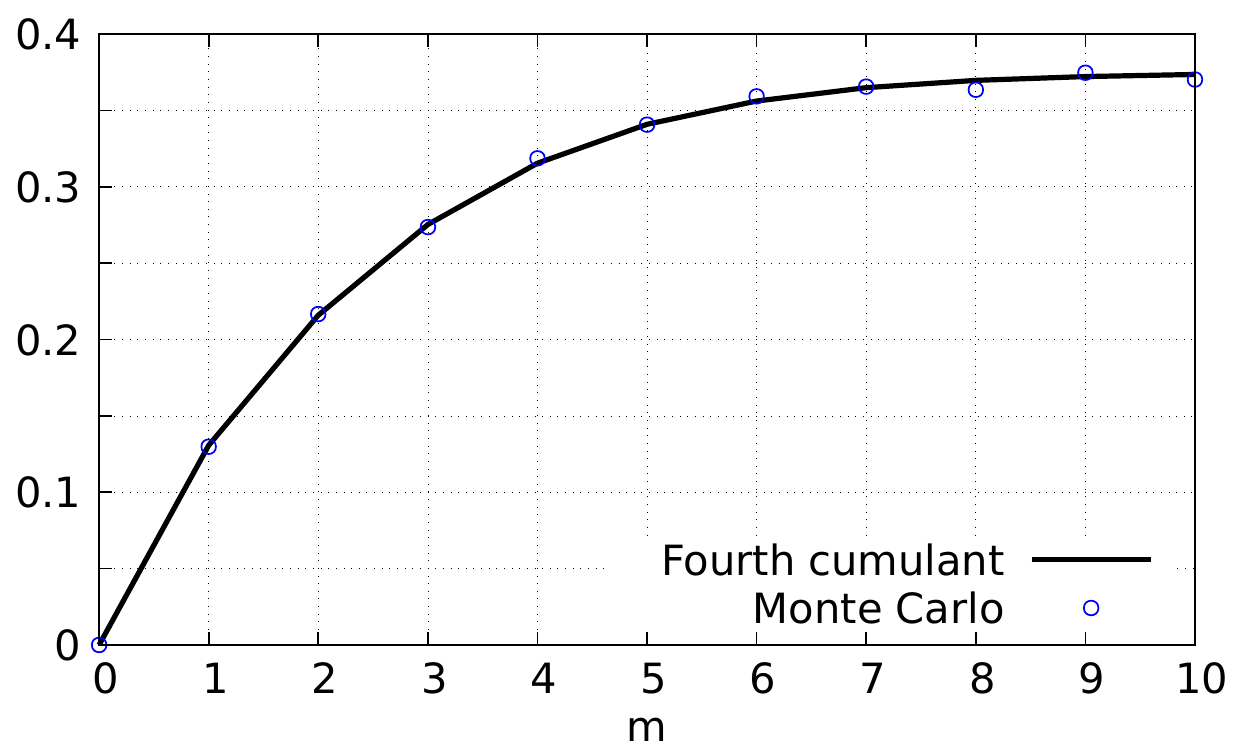}
    \caption{Fourth cumulant of $X(m)$.} 
 \end{subfigure}
 \begin{subfigure}[b]{0.49\textwidth}
    \includegraphics[width=1\linewidth, height=5cm]{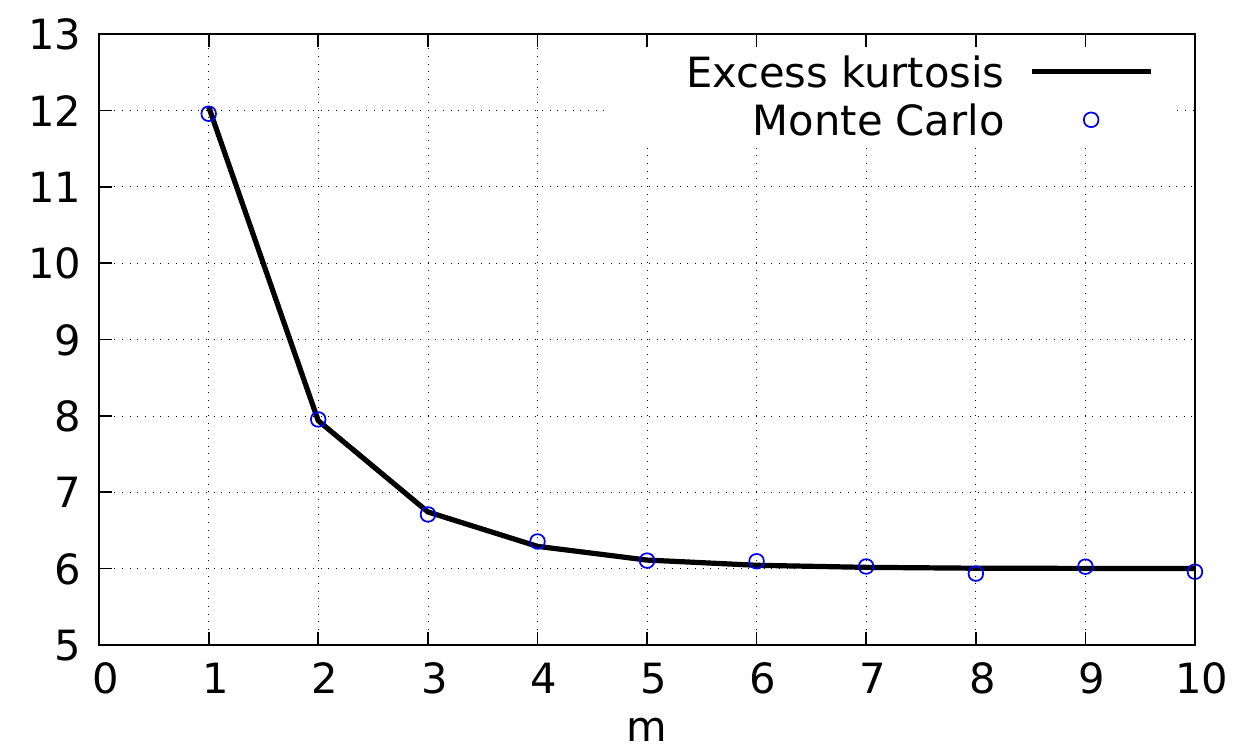} 
 \caption{Excess kurtosis $\kappa^{(4)} (m)/\big( \kappa^{(2)} (m)\big)^2$.} 
 \end{subfigure}
 \caption{Fourth cumulant and kurtosis of $X(m)$.} 
 \label{fig7} 
\end{figure}

 \vspace{-0.4cm} 
 \noindent 
 Higher order moments of $X(m)$ 
 can be obtained from Corollary~\ref{p6.5}
 using Mathematica, and are plotted
 with $f_{l,m}(s)=s$ at the orders $3$ and $4$ in Figures~\ref{fig6}-\ref{fig7}, 
 along with Monte Carlo simulations used for
 confirmation. 

\footnotesize 

\def\cprime{$'$} \def\polhk#1{\setbox0=\hbox{#1}{\ooalign{\hidewidth
  \lower1.5ex\hbox{`}\hidewidth\crcr\unhbox0}}}
  \def\polhk#1{\setbox0=\hbox{#1}{\ooalign{\hidewidth
  \lower1.5ex\hbox{`}\hidewidth\crcr\unhbox0}}} \def\cprime{$'$}

\end{document}